\documentclass[12pt]{amsart}
\usepackage{amsfonts,amscd,amsmath,amssymb,amsthm,mathrsfs}
\usepackage{array,latexsym,color,float,enumerate}
\usepackage{graphicx,epsfig}
\usepackage[top=2.7cm, bottom=2.7cm, left=3.2cm, right=3.2cm]{geometry}
\usepackage[all]{xy}
\usepackage[hypertex,linkbordercolor={0 0 1}]{hyperref}

\makeindex

\usepackage{marginnote}

\theoremstyle{plain}
\newtheorem{theorem}{Theorem}[section]
\newtheorem{lemma}[theorem]{Lemma}
\newtheorem{proposition}[theorem]{Proposition}
\newtheorem{corrolary}[theorem]{Corollary}
\theoremstyle{definition}
\newtheorem{definition}[theorem]{Definition}
\newtheorem{example}[theorem]{Example}
\newtheorem{noname}[theorem]{}
\newtheorem{remark}[theorem]{Remark}
\newtheorem{construction}[theorem]{Construction}
\newtheorem{notation}[theorem]{Notation}
\theoremstyle{remark}
\newtheorem*{smallremark}{Remark}

\newtheorem{case}{Case} \makeatletter \@addtoreset{case}{theorem}\makeatother
\newtheorem{claim}{Claim} \makeatletter \@addtoreset{claim}{theorem}\makeatother

\newcommand{\bthm}{\begin{theorem}}
\newcommand{\bprop}{\begin{proposition}}
\newcommand{\blem}{\begin{lemma}}
\newcommand{\bcor}{\begin{corrolary}}
\newcommand{\brem}{\begin{remark}}
\newcommand{\bdfn}{\begin{definition}}
\newcommand{\bitem}{\begin{itemize}}
\newcommand{\bex}{\begin{example}}
\newcommand{\bno}{\begin{noname}}
\newcommand{\bsrem}{\begin{smallremark}}
\newcommand{\bnot}{\begin{notation}}
\newcommand{\bcon}{\begin{construction}}
\newcommand{\bca}{\begin{case}}
\newcommand{\bcl}{\begin{claim}}

\newcommand{\ecl}{\end{claim}}
\newcommand{\eca}{\end{case}}
\newcommand{\econ}{\end{construction}}
\newcommand{\enot}{\end{notation}}
\newcommand{\esrem}{\end{smallremark}}
\newcommand{\eno}{\end{noname}}
\newcommand{\eex}{\end{example}}
\newcommand{\eitem}{\end{itemize}}
\newcommand{\ethm}{\end{theorem}}
\newcommand{\eprop}{\end{proposition}}
\newcommand{\elem}{\end{lemma}}
\newcommand{\ecor}{\end{corrolary}}
\newcommand{\erem}{\end{remark}}
\newcommand{\edfn}{\end{definition}}
\newcommand{\benum}{\begin{enumerate}}
\newcommand{\eenum}{\end{enumerate}}
\newcommand{\un}{\underline}
\newcommand{\ov}{\overline}

\newcommand{\wt}{\widetilde}
\newcommand{\cal}[1]{\mathcal{#1}}

\def\8{\infty}
\def\.{\cdot}
\def\PP{\mathbb{P}}
\def\F{\mathbb{F}}
\def\C{\mathbb{C}}
\def\Z{\mathbb{Z}}
\def\N{\mathbb{N}}
\def\Q{\mathbb{Q}}

\def\E{\widehat{E}}
\def\ovk{\overline\kappa}

\def\:{\colon}
\def\map{\dashrightarrow}
\def\transf{\rightsquigarrow}

\def\Coker{\operatorname{Coker}}

\def\Supp{\operatorname{Supp}}
\def\Pic{\operatorname{Pic}}
\def\dim{\operatorname{dim}}
\def\Exc{\operatorname{Exc}}

\begin{document}
\title[Singular $\Q$-homology planes II]{Classification of singular $\Q$-homology planes. \\ II. $\C^1$- and $\C^*$-rulings.}
\author[Karol Palka]{Karol Palka}
\address{Karol Palka: Institute of Mathematics, University of Warsaw, ul. Banacha 2, 02-097 Warsaw, Poland}
\address{Institute of Mathematics, Polish Academy of Sciences, ul. \'{S}niadeckich 8, 00-956 Warsaw, Poland}
\thanks{The author was supported by Polish Grant NCN 6086/B/H03/2011/40}
\email{palka@impan.pl}

\begin{abstract} A $\Q$-homology plane is a normal complex algebraic surface having trivial rational homology. We classify singular $\Q$-homology planes which are $\C^1$- or $\C^*$-ruled. We analyze their completions, the number of different rulings, the number of affine lines on it and we give constructions. Together with previously known results this completes the classification of $\Q$-homology planes with smooth locus of non-general type. We show also that the dimension of a family of homeomorphic but non-isomorphic singular $\Q$-homology planes having the same weighted boundary, singularities and Kodaira dimension can be arbitrarily big. \end{abstract}

\maketitle

We work with complex algebraic varieties.

\section{Main results}

A $\Q$-homology plane is a normal surface whose rational cohomology is the same as that of $\C^2$. The paper is the last piece of the classification of $\Q$-homology planes having smooth locus of non-general type. The classification is built on work of many authors, for summary of what is known for smooth and singular $\Q$-homology planes see \cite[\S 3.4]{Miyan-OpenSurf} and \cite{Palka-recent_progress}. In \cite{Palka-classification1} we have classified singular $\Q$-homology planes with non-quotient singularities showing in particular that they are quotients of affine cones over projective curves by actions of finite groups which respect the set of lines through the vertex. In \cite{Palka-k(S_0)=0} we classified singular $\Q$-homology planes whose smooth locus is of non-general type and admits no $\C^1$- and no $\C^*$-ruling (so-called \emph{exceptional planes}). Here we classify singular $\Q$-homology planes which admit a $\C^1$- or a $\C^*$-ruling. We analyze completions and boundaries rather than the open surfaces themselves. To deal with non-uniqueness of these we use the notion of a \emph{balanced} and a \emph{strongly balanced} weighted boundary and completion of an open surface (cf. \ref{def2:balanced}, \ref{def2:normalization of boundary}).

The classification of $\C^1$- and $\C^*$-ruled $\Q$-homology planes goes by giving necessary and sufficient conditions for a $\C^1$- or $\C^*$-ruled open surface to be a $\Q$-homology plane (see \ref{lem5:when ruling gives sQhp}, \ref{lem5:when d(D)=0} and the remarks before) and then giving a general construction (see section \ref{ssec:construction}). We compute the Kodaira dimension of a $\C^*$-ruled singular $\Q$-homology plane and of its smooth locus (\ref{thm5:(k,k0) computations}) in terms of properties of singular fibers and then we list the planes with smooth locus of Kodaira dimension zero (section \ref{ssec:k(S_0)=0}). As a corollary from the classification we obtain in particular the following result.

\bthm\label{thm5:singularities} Let $S'$ be a singular $\Q$-homology plane and let $S_0$ be its smooth locus. Assume $S'$ is not affine-ruled and $\ovk(S_0)\neq 2$. Then the following hold. \benum[(1)]

\item Either $S'$ has a unique balanced completion up to isomorphism or it admits an untwisted $\C^*$-ruling with base $\C^1$ and more than one singular fiber. In the last case $S'$ has exactly two strongly balanced completions.

\item If $S'$ has more than one singular point then it has exactly two singular points, both of Dynkin type $A_1$, and there is a twisted $\C^*$-ruling of $S'$, such that both singular points are contained in a unique fiber isomorphic to $\C^1$.

\item If $S'$ contains a quotient non-cyclic singularity then either $S'\cong\C^2/G$ for a small finite noncyclic subgroup of $GL(2,\C)$ or $S'$ has a twisted $\C^*$-ruling. In the last case the unique fiber isomorphic to $\C^1$ is of type (A)(iv) (cf. \ref{thm5:(k,k0) computations}) and contains a singular point of Dynkin type $D_k$ for some $k\geq 4$. \eenum \ethm

We now comment on other corollaries from the obtained classification. First, the case when $S'$ has exactly one singular point and it is a cyclic singularity can occur. Second, we show that if $S'$ is affine-ruled then its strongly balanced weighted boundary is unique unless it is a chain, but even if it is unique there still may be infinitely many strongly balanced completions (see  \ref{ex:one boundary many completions}). Third, the singularities of affine-ruled $S'$ are necessarily cyclic but there may be arbitrarily many of them (see \cite{MiSu-hPlanes} or section \ref{sec:k(S0) negative}). As for the remaining case $\ovk(S_0)=2$ let us mention that it follows from the logarithmic Bogomolov-Miyaoka-Yau inequality (cf. \cite{Palka-classification1}) that $S'$ has only one singular point and it is of quotient type.

It is known (cf. \cite{FZ-deformations}) that smooth $\Q$-homology planes can have moduli. The same happens for singular ones. We prove the following result.

\bthm\label{thm:moduli} There exist arbitrarily high-dimensional families of non-isomorphic singular $\Q$-homology planes having negative Kodaira dimension of the smooth locus, the same singularities, homeomorphism type and the same weighted strongly balanced boundary. \ethm

It is an important property of any $\Q$-homology plane with smooth locus of general type that it does not contain topologically contractible curves. In fact the number of contractible curves on a $\Q$-homology plane is known except in the case when the surface is singular and the smooth locus has Kodaira dimension zero (see the discussion in section \ref{ssec:number of C*-rulings}). In theorem \ref{thm5:uniqueness of C*-rulings} we compute the number of different $\C^*$-rulings a $\Q$-homology plane can have. The computation of the number of contractible curves follows from it.

\bthm\label{cor5:number of contractible curves} If a singular $\Q$-homology plane has smooth locus of Kodaira dimension zero then it contains one or two irreducible topologically contractible curves in case the smooth locus admits a $\C^*$-ruling and contains no such curves otherwise.\ethm

We now comment on some related results of other authors. First of all the notion of a balanced weighted boundary and a balanced completion of an open surface  (cf. \ref{def2:normalization of boundary}) is a more flexible version of the notion of a 'standard graph' from \cite{FKZ-weighted-graphs}, which has its origin in the paper of \cite{Daigle-weighted-graphs}. It follows from above that every $\Q$-homology plane admits up to isomorphism one or two strongly balanced boundaries, it is not so for the standard ones. The set of such boundaries is a useful invariant of the surface. Second, integral homology groups and necessary conditions for singular fibers of $\C^1$- and $\C^*$-ruled $\Q$-homology planes have been already analyzed in \cite{MiSu-hPlanes}. However, as we show, for $\C^*$-rulings these conditions are not sufficient (cf. \ref{ex: C*-ruled with d(D)=0:I}, \ref{ex: C*-ruled with d(D)=0:II}) and a more detailed analysis is necessary. Moreover, formulas for the Kodaira dimension in terms of singular fibers from loc. cit. require nontrivial corrections (see section \ref{ssec:Kodaira dimensions}).

\section{Preliminaries}\label{sec:preliminaries}

We follow the notational conventions and terminology of \cite{Miyan-OpenSurf}, \cite{Fujita} and \cite{Palka-classification1}. We recall some of them for the convenience of the reader.

\subsection{Divisors and normal pairs}\label{ssec:divs and pairs} Let $T=\sum t_iT_i$ be an snc-divisor on a smooth complete surface with distinct irreducible components $T_i$. Then $\un T=\sum T_i$, where the sum runs over $i$ with $t_i\neq 0$, is the reduced divisor with the same support as $T$ and $\beta_T(T_i)=\un T\cdot (\un T-T_i)$ is called the \emph{branching number} of $T_i$. A \emph{tip} has $\beta_T(T_i)\leq 1$. By $Q(T)$ we denote the intersection matrix of $T$, we put $d(0)=1$ and $d(T)=det(-Q(T))$ for $T\neq 0$. Numerical equivalence of divisors is denoted by $\equiv.$

If $T$ is reduced and its dual graph is linear, it is called a \emph{chain} and writing it as a sum of irreducible components $T=T_1+\ldots+T_n$ we assume that $T_i\cdot T_{i+1}=1$ for $1\leq i\leq n-1$. We put $T^t=T_n+\ldots+T_1$. If $T$ is a rational chain then we write $T=[-T_1^2,\ldots,-T_n^2]$. A rational chain with all $T_i^2\leq -2$ is called \emph{admissible}. A \emph{fork} (\emph{comb}) is a rational tree with a unique branching component, the component has $\beta=3$ ($\beta\geq 3$).

Let $D$ be some reduced snc-divisor, which is not an admissible chain. A rational chain with support contained in $D$, not containing branching components of $D$ and containing one of its tips is called a \emph{twig} of $D$. For an admissible (ordered) chain we put $$e(T)=\frac{d(T-T_1)}{d(T)}\text{\ and\ }\wt e(T)=e(T^t).$$ In general $e(T)$ and $\wt e(T)$ are defined as the sums of respective numbers computed for all maximal admissible twigs of $T$. Here the convention that the tip of the twig is the first component is used.

If $X$ is a complete surface and $D$ a reduced snc-divisor contained in the smooth part of $X$ then we call $(X,D)$ an \emph{snc-pair} and we write $X-D$ for $X\setminus D$. The pair is \emph{normal} (\emph{smooth}) if $X$ is normal (resp. smooth). If $X$ is a normal surface then an embedding $\iota\:X\to \ov X$, where $(\ov X,\ov X\setminus X)$ is a normal pair, is called a \emph{normal completion} of $X$. If $X$ is smooth then $\ov X$ is smooth and $(\ov X,D,\iota)$ is called a \emph{smooth completion} of $X$. A morphism of two completions $\iota_j:X\to \ov X_j$, $j=1,2$ of a given surface $X$ is a morphism $f\:\ov X_1\to \ov X_2$, such that $\iota_2=f\circ\iota_1$.

Let $\pi\:(X,D)\to (X',D')$ be a birational morphism of normal pairs. We put $\pi^{-1}D'=\un {\pi^*D'}$, i.e. $\pi^{-1}D'$ is the reduced total transform of $D'$. Assume $\pi^{-1}D'=D$. If $\pi$ is a blow-up then we call it \emph{subdivisional} (\emph{sprouting}) for $D'$ if its center belongs to two (one) components of $D'$. In general we say that $\pi$ is \emph{subdivisional} for $D'$ (and for $D$) if for any component $T$ of $D'$ we have $\beta_{D'}(T)=\beta_D(\pi^{-1}T)$. The exceptional locus of a birational morphism between two surfaces $\eta:X\to X'$, denoted by $\Exc(\eta)$, is defined as the locus of points in $X$ for which $\eta$ is not a local isomorphism.

A $b$-curve is a smooth rational curve with self-intersection $b$. A divisor is snc-minimal if all its $(-1)$-curves are branching. We write $K_X$ for the canonical divisor on a complete surface $X$.

\bdfn\label{def2:connected morphism} A birational morphism of surfaces $\pi:X\to X'$ is a \emph{connected modification} if it is proper, $\pi(\Exc(\pi))$ is a smooth point on $X'$ and $\Exc(\pi)$ contains a unique $(-1)$-curve. In case $\pi$ is a morphism of pairs $\pi\:(X,D)\to (X',D')$, such that $\pi^{-1}(D')=D$ and $\pi(\Exc(\pi))\in D'$, then we call it a \emph{connected modification over $D'$}.\edfn

A sequence of blow-downs (and its reversing sequence of blow-ups) whose composition is a connected modification will be called a \emph{connected sequence of blow-downs (blow-ups)}.

\subsection{Rational rulings}\label{ssec:rulings}

A surjective morphism $p_0:X_0\to B_0$ of a normal surface onto a smooth curve is a \emph{rational ruling} if general fibers are rational curves. By a \emph{completion of $p_0$} we mean a triple $(X,D,p)$, where $(X,D)$ is a normal completion of $X_0$ and $p\:X\to B$ is an extension of $p_0$ to a $\PP^1$-ruling with $B$ being a smooth completion of $B_0$. We say that $p$ is a \emph{minimal completion of $p_0$} if $p$ does not dominate any other completion of $p_0$. In this case we also say that $D$ is \emph{$p$-minimal}. It is easy to check that $D$ is $p$-minimal if and only if all its non-branching $(-1)$-curves are horizontal. Let $F$ be a fiber of $p$. An irreducible curve $G\subseteq X$ is an \emph{$n$-section} of $p$ if $G\cdot F=n$. A \emph{section} is a $1$-section. We call $p_0$ a $\C^{(n*)}$-ruling if $F\cdot D=n+1$, $n\geq 1$. In case $n=0$ we call it a $\C^1$-ruling or an \emph{affine ruling}. The arithmetic genus of $F$ ($p_a(F)=\frac{1}{2}F\cdot (K_X+F)$) vanishes and $F^2=0$. Conversely, it is well-known that an effective divisor with these properties on a complete surface is a fiber of such a ruling (see \cite[V.4.3]{BHPV}). If $J$ is a component of $F$ then we denote by $\mu_F(J)$ the multiplicity of $J$, i.e. $F=\mu_F(J)J+F'$, where $F'$ is effective and $J\not\subseteq F'$. The structure of fibers of a $\PP^1$-ruling is well known (see \cite[\S 4]{Fujita}).

\blem\label{lem2:fiber with one exc comp} Let $F$ be a singular fiber of a $\PP^1$-ruling of a smooth complete surface. Then $F$ is a tree of rational curves and it contains a $(-1)$-curve. Each $(-1)$-curve of $F$ meets at most two other components. If $F$ contains a unique $(-1)$-curve $C$ then:\benum[(i)]

\item $\mu(C)>1$ and there are exactly two components of $F$ with multiplicity one, they are tips of the fiber,

\item if $\mu(C)=2$ then either $F=[2,1,2]$ or $C$ is a tip of $F$ and then $\un F-C=[2,2,2]$ or  $\un F-C$ is a $(-2)$-fork of type $(2,2,n)$,

\item if $\un F$ is not a chain then the connected component of $\un F-C$ not containing curves of multiplicity one is a chain (possibly empty).\eenum\elem

We define $$\Sigma_{X-D}=\underset{F\nsubseteq D}{\sum}(\sigma(F)-1),$$ where $\sigma(F)$ is the number of $(X-D)$-components of a fiber $F$ (cf. \cite[4.16]{Fujita}). If $p$ is a $\PP^1$-ruling as above then we say that an irreducible curve $G$ is \emph{vertical (for $p$)} if $p_*G=0$, otherwise it is \emph{horizontal}. A divisor is vertical (horizontal) if all its components are vertical (horizontal). We decompose $D$ as $D=D_h+(D-D_h)$, where $D_h$ is horizontal and $D-D_h$ is vertical. The numbers $h$ and $\nu$ are defined respectively as the number of irreducible components of $D_h$ and as the number of fibers contained in $D$. We have (cf. \S 4 loc. cit.): $$\Sigma_{X-D}=h+\nu+b_2(X)-b_2(D)-2.$$ We call a connected component of $F\cap D$ a \emph{$D$-rivet} (or \emph{rivet} if this makes no confusion) if it meets $D_h$ at more than one point or if it is a node of $D_h$.

\bdfn\label{def2:(un)twisted and columnar fiber} Let $(X,D,p)$ be a completion of a $\C^*$-ruling of a normal surface $X$. We say that the original ruling $p_0=p_{|X-D}$ is \emph{twisted} if $D_h$ is a $2$-section. If $D_h$ consists of two sections we say that $p_0$ is \emph{untwisted}. Let $F$ be a singular fiber $F$ of $p$ which does not contain singular points of $X$. We say that $F$ is \emph{columnar} if and only if $\un F$ is a chain which can be written as $$\un F=A_n+\ldots+A_1+C+B_1+\ldots+B_m,$$ where $C$ is a unique $(-1)$-curve and $D_h$ meets $F$ exactly in $A_n$ and $B_m$. The chains $A=A_1+\ldots+A_n$ and $B=B_1+\ldots+B_m$ are called \emph{adjoint chains}. \edfn

\bsrem By expansion properties of determinants (cf. \cite[2.1.1]{KR-ContrSurf}) and the fact that $d(A)$ and $d(A-A_1)$ are coprime we have $e(A)+e(B)=1$ and $d(A)=d(B)=\mu_F(C)$. In fact we have also $\wt e(B)+\wt e(A)=1$ (see \cite[3.7]{Fujita}). \esrem

\subsection{Balanced completions}\label{ssec:completions and boundaries}

\bdfn\label{def2:weighted boundary} A pair $(D,w)$ consisting of a complete curve $D$ and a rationally-valued function $w$ defined on the set of irreducible components of $D$ is called a \emph{weighted curve}. If $(X,D)$ is a normal pair then $(D,w)$ with $w$ defined by $w(D_i)=D_i^2$ is a \emph{weighted boundary} of $X-D$.\edfn

\bdfn\label{def2:flow} Let $(X,D)$ be a normal pair.\benum [(i)]

\item Let $L$ be a $0$-curve which is a non-branching component of $D$ and let $c\in L$ be chosen so that if $L$ intersects two other components of $D$ then $c$ is one of the points of intersection. Make a blow-up of $c$ and contract the proper transform of $L$. The resulting pair $(X',D')$, where $D'$ is the reduced direct image of the total transform of $D$ is called an \emph{elementary transform of $(X,D)$}. The pair $\Phi=(\Phi^\circ,\Phi^\bullet)$ consisting of an assignment $\Phi^\circ\:(X,D)\mapsto (X',D')$ together with the resulting rational mapping $\Phi^\bullet\:X\map X'$ is called an \emph{elementary transformation over $D$}. $\Phi$ is \emph{inner (for $D$)} if $\beta_D(L)=2$ and \emph{outer (for $D$)} if $\beta_D(L)=1$. The point $c\in L$ is the \emph{center} of $\Phi$.

\item For a sequence of (inner) elementary transformations $\Phi^\circ_i\:(X_i,D_i)\mapsto (X_{i+1},D_{i+1})$, $i=1,\ldots,{n-1}$ we put $\Phi^\circ=(\Phi^\circ_1,\ldots,\Phi^\circ_{n-1})$, $\Phi^\bullet=(\Phi^\bullet_1,\ldots,\Phi^\bullet_{n-1})$ and we call $\Phi=(\Phi^\circ,\Phi^\bullet)$ \emph{an (inner) flow in $D_1$}. We denote it by $\Phi\:(X_1,D_1)\transf (X_n,D_n)$. \eenum\edfn

Note that $\Phi^\bullet=(\Phi^\bullet_1,\ldots,\Phi^\bullet_{n-1})$ induces a rational mapping $X_1\map X_n$, which we also denote by $\Phi^\bullet$. There exists the largest open subset of $X_1$ on which $\Phi_1^\bullet$ is a morphism, the complement of this subset is called the \emph{support of $\Phi$}. Clearly, $\Supp \Phi_1\subseteq D_1$. If $\Supp \Phi=\emptyset$ then $\Phi$ is a \emph{trivial flow}.

A weighted curve $(D,w)$ determines the weighted dual graph of $D$. If $(D,w)$ is a weighted boundary coming from a fixed normal pair $(X,D)$ we omit the weight function $w$ from the notation. Note that for $\Phi$ as above $D_1$ and $D_n$ are isomorphic as curves. They have the same dual graphs, but usually different weights of components.

\bex\label{ex2:flow} Let $T=[0,0,a_1,\ldots,a_n]$. Then each chain of type $[0,b,a_1,\ldots,a_n]$, $[a_1,\ldots,a_{k-1},a_k-b,0,b,a_{k+1},\ldots,a_n]$ or $[a_1,\ldots,a_n,b,0]$ where $1\leq k\leq n$ and $b\in\Z$, can be obtained from $T$ by a flow. This follows easily from the observation that an elementary transformation changes the chains $[w,x,0,y-1,z]$ and $[w,x-1,0,y,z]$ one into another. Looking at the dual graph we see the weights can 'flow' from one side of a $0$-curve to another, including the possibility that they vanish ($b=0$ or $b=a_k$). If they do then again the weights can flow through the new zero. \eex

\bdfn\label{def2:balanced} A rational chain $D=[a_1,\ldots,a_n]$ is \emph{balanced} if $a_1,\ldots,a_n\in\{0,2,3,\ldots\}$ or if $D=[1]$. A reduced snc-divisor whose dual graph contains no loops (snc-forest) is \emph{balanced} if all rational chains contained in $D$ which do not contain branching components of the divisor are balanced. A normal pair $(X,D)$ is \emph{balanced} if $D$ is balanced. \edfn

Recall that if $(X_i,D_i)$ for $i=1,2$ are normal pairs such that $X_1-D_1\cong X_2-D_2$ then $D_1$ is a forest if and only if $D_2$ is a forest.

\bprop\label{prop2:boundary uniqueness} Any normal surface which admits a normal completion with a forest as a boundary has a balanced completion. Two such completions differ by a flow. \eprop

As we discovered after completing the proof, the above proposition in a more general version was proved in a graph theoretic context in \cite{FKZ-weighted-graphs} (see Theorem 3.1 and Corollary 3.36 loc. cit.). We leave therefore our more direct arguments to be published elsewhere. In fact, some key observations were done earlier in \cite{Daigle-weighted-graphs} (see 4.23.1, 3.2, 5.2 loc. cit.). Let us restate some definitions from \cite{FKZ-weighted-graphs} on the level of pairs.

\bdfn\label{def2:standard boundary} Let $(X,D)$ be a normal pair and assume $D$ is an snc-forest.\benum[(i)]

\item Connected components of the divisor which remains after subtracting all non-rational and all branching components of $D$ are called the \emph{segments of $D$}.

\item $D$ is \emph{standard} if for each of its connected components either this component is equal to $[1]$ or all its segments are of types $[0]$, $[0,0,0]$ or $[0^{2k},a_1,\ldots,a_n]$ with $k\in \{0,1\}$ and $a_1,\ldots,a_n\geq 2$.

\item if $D_0=[0,0,a_1,\ldots,a_n]$ with $a_i\neq 0$, $i=1,\ldots,n$ is a segment of $D$ then a \emph{reversion of $D_0$} is a nontrivial flow $\Phi:(X,D)\transf (X',D')$ with support in $D_0$, which is inner for $D_0$ and for which $D'-(\Phi^\bullet)_*(D-D_0)=[a_1,a_2,\ldots,a_n,0,0]$.
\eenum \edfn

The condition that $\Phi$ is nontrivial is introduced for the following reason: we want the reversion to transform the two zeros 'to the other end' of the chain, and the condition in necessary to force this in case $D$ is symmetric, i.e when $[a_1,\ldots,a_n]^t=[a_1,\ldots,a_n]$. Standard chains are called \emph{canonical} in \cite{Daigle-weighted-graphs}. Note that the Hodge index theorem implies that if $(X,D)$ is a smooth pair and $D$ is a forest then it cannot have segments of type $[0^{2k+1}]$ or $[0^{2k},a_1,\ldots,a_n]$ for $k>1$ and can have at most one such segment with $k=1$.

Clearly, not every balanced forest is standard, but by a flow one can easily change it to such. Now it follows from \ref{prop2:boundary uniqueness} that if $D$ and $D'$ are two standard boundaries of the same surface and $D$ is a chain then either $D$ and $D'$ are isomorphic as weighted curves or $D'$ is the reversion of $D$. Unfortunately, the notion of a standard boundary in not as restrictive as one may imagine and the difference between two standard boundaries can be more than just a reversion of some segments. An additional ambiguity is related to the existence of segments of type $[0^{2k+1}]$. Namely, if $[0^{2k+1}]$ is a segment of $D$ then one can change by a flow the self-intersections of the components of $D$ intersecting the segment. For example, consider a surface whose standard boundary is a rational fork with a dual graph $$ \xymatrix{{-2}\ar@{-}[r] &{b}\ar@{-}[r]\ar@{-}[d]& {-2}\\ {} & {0} & {} } $$ for some $b\in\Z$. Then for any $b\in \Z$ there is a completion of this surface for which the boundary is standard and has the dual graph as above.\footnote[1]{This observation was missed in \cite{FKZ-weighted-graphs} and the corollary 3.33 loc. cit. is false. See \cite{FKZ-correction} for corrections. In \cite[Solution to problem 5, p. 45]{Daigle-weighted-graphs} this ambiguity is implicitly taken into account without restricting to balanced divisors. } We therefore introduce the following more restrictive conditions.

\bdfn\label{def2:normalization of boundary} A balanced snc-forest $D$ is \emph{strongly balanced} if and only if it is standard and either $D$ contains no segments of type $[0]$, $[0,0,0]$ or for at least one of such segments there is a component $B\subseteq D$ intersecting it, such that $B^2=0$. A normal pair $(X,D)$ for which $D$ is a forest is \emph{strongly balanced} if $D$ is strongly balanced.\edfn

\subsection{Basic properties of $\Q$-homology planes}\label{ssec:Qhp notation}

We assume that $S'$ is a \emph{singular $\Q$-homology plane}, i.e. a normal non-smooth complex algebraic surface with $H^*(S',\Q)\cong \Q$. Let $\epsilon \colon S \to S'$ be a resolution such that the inverse image of the singular locus is an snc-divisor and let $(\ov S,D)$ be a smooth completion of $S$. Denote the singular points of $S'$ by $p_1,\ldots,p_q$ and the smooth locus by $S_0$. We put $\E_i=\epsilon^{-1}(p_i)$ and we assume that $\E =\E_1+\E_2+\ldots +\E_q$ is snc-minimal. Recall that $S'$ is called \emph{logarithmic} if and only if every singular point of $S'$ is locally analytically isomorphic to $\C^2/G$ for some subgroup $G<GL(2,\C)$ (a 'quotient' singularity). In \cite{Palka-classification1} we classified non-logarithmic $\Q$-homology planes. In particular it is known that they do not admit $\C^1$- or $\C^*$-rulings. Therefore, from now on we assume that $S'$ is logarithmic. It follows that each $\E_i$ is either an admissible chain or an admissible fork (i.e. an snc-minimal fork with negative definite intersection matrix). By \cite{GPS-logQhp_rational} $S'$ is rational. By the argument in \cite[2.4]{Fujita} it is affine.

\bprop\label{prop:homology} Let the notation be as above. Then:\benum[(i)]

\item $D$ is a rational tree with $d(D)=-d(\E)\cdot|H_1(S',\Z)|^2$,

\item the embedding $D\cup \E\to \ov S$ induces an isomorphism on $H_2(-,\Q)$,

\item $\pi_1(S')\cong \pi_1(S)$ and $H_k(S',\Z)=0$ for $k>1$,

\item $b_i(S_0)=0$ for $i=1,2,4$, $b_3(S_0)=q$,

\item $\Sigma_{S_0}=h+\nu-2$ and $\nu\leq 1$. \eenum\eprop

\begin{proof} See \cite[3.1, 3.2]{Palka-classification1} and \cite[2.2]{MiSu-hPlanes}.  \end{proof}

We have the following criterion.

\blem\label{lem5:when ruling gives sQhp} Let $(\ov S,T)$ be a smooth pair and let $p\:\ov S\to \PP^1$ be a $\PP^1$-ruling. Assume the following conditions are satisfied:\benum[(i)]

\item there exists a unique connected component $D$ of $T$ which is not vertical,

\item $D$ is a rational tree,

\item $\Sigma_{\ov S-T}=h+\nu-2$,

\item $d(D)\neq 0$.

\eenum Then the surface $S'$ defined as the image of $\ov S-D$ after contraction of connected components of $T-D$ to points is a rational $\Q$-homology plane and $p$ induces a rational ruling of $S'$. Conversely, if $p'\:S'\to B$ is a rational ruling of a rational $\Q$-homology plane $S'$ then any completion $(\ov S,T,p)$ of the restriction of $p'$ to the smooth locus of $S'$ has the above properties.\elem

\begin{proof} Since the base of $p$ has some component of $D$ as a branched cover, it is rational, hence $\ov S$ is rational. We may assume $T$ is $p$-minimal. Put $\E=T-D$. Since $\E$ is vertical and since $\E\cap D=\emptyset$, $Q(\E)$ is negative definite and $b_1(\E)=0$. Fujita's equation $$\Sigma_{\ov S-T}=h+\nu-2+b_2(\ov S)-b_2(D+\E)$$ gives $b_2(\ov S)=b_2(T)$, so by (iv) the inclusion $T\to \ov S$ induces an isomorphism on $H_2(-,\Q)$. By \cite[2.4]{Palka-classification1} $S'$ is normal and affine, in particular $b_4(S')=b_3(S')=0$. Since $b_1(D)=0$, the exact sequence of the pair $(\ov S,D)$ together with the Lefschetz duality give $$b_2(S)=b_2(\ov S,D)=b_2(\ov S)-b_2(D)=b_2(\E).$$ Since $b_1(\E)=0$, we get from the exact sequence of the pair $(S,\E)$ that $b_2(S')=b_2(S,\E)=b_2(S)-b_2(\E)=0$. Now $$\chi(S')=\chi(\ov S)-\chi(D\cup\E)+b_0(\E)=b_0(D)=1,$$ so we obtain $b_1(S')=b_2(S')=0$, hence $S'$ is $\Q$-acyclic.

Conversely, if $p'$ is as above then let $\E$ be an exceptional divisor of a resolution of singularities of $S'$ and let $D=T-\E$. Since $\E$ is vertical for the $\PP^1$-ruling $p$, we have $b_1(\E)=0$. Then the necessity of the above conditions follows from 3.1. and 3.2 loc. cit. \end{proof}

\section{\label{sec:k(S0) negative} Smooth locus of negative Kodaira dimension}

In this section we assume that the smooth locus $S_0$ of the logarithmic $\Q$-homology plane $S'$ has negative Kodaira dimension. This implies that the Kodaira dimension of $S'$ is also negative. The case was analyzed in \cite[2.5-2.8]{MiSu-hPlanes}, where a structure theorem was given. We, in particular, recover these results in \ref{lem4:affine-ruled S'} and \ref{prop4:non affine-ruled S'}, but we concentrate on analyzing possible completions and boundaries instead of $S'$ itself. This gives more information, allows to give a construction and to answer the question of uniqueness of an affine ruling of $S_0$ (in case it exists). The information about completions is also used in the analysis of an example where moduli occur.

\bprop\label{prop4:non affine-ruled S'} If a singular $\Q$-homology plane has smooth locus of negative Kodaira dimension then it is affine-ruled or isomorphic to $\C^2/G$ for some small, noncyclic subgroup $G<GL(2,\C)$. The surfaces $\C^2/G$ and $\C^2/G'$ are isomorphic if and only if $G$ and $G'$ are conjugate in $GL(2,\C)$. The minimal normal completion of $\C^2/G$ is unique, the boundary is a non-admissible rational fork with admissible twigs. \eprop

\begin{proof} For the first part of the statement we follow the arguments of \cite[\S 3]{KR-ContrSurf}. Assume that $S'$ is not affine-ruled. Then $S_0$ is not affine-ruled. Since $S'$ is affine, $D+\E$ is not negative definite, so by \cite[2.5.1]{Miyan-OpenSurf} $S_0$ contains a Platonically $\C^*$-fibred open subset $U$, which is its almost minimal model. Moreover, $\chi(U)\leq\chi(S_0)$ (cf. \cite[2.8]{Palka-k(S_0)=0}). The algorithm of construction of an almost minimal model (see \cite[2.3.8, 2.3.11]{Miyan-OpenSurf}) implies that $S_0-U$ is a disjoint sum of $s$ curves isomorphic to $\C$ and $s'$ curves isomorphic to $\C^*$ for some $s,s'\in\N$. It follows that $$0=\chi(U)=\chi(S_0)-s=\chi(S')-q-s=1-q-s,$$ so $s=0$, $q=1$ and $s'\leq1$. If $s'\neq 0$ then the boundary divisor of $U$ is connected, hence $U$ and $S_0$ are affine-ruled. Thus $s'=0$, $S_0=U$ and by \cite{MiTs-logDelPezzo} $S'\cong\C^2/G$, where $G$ is a small noncyclic subgroup of $GL(2,\C)$.

Suppose $G$ and $G'$ are two subgroups of $GL(2,\C)$, such that $\C^2/G\cong\C^2/G'$. Then $\widehat{\cal O}_{\C^2/G,(0)}\cong \widehat{\cal O}_{\C^2/G',(0)}$, so if $G$ and $G'$ are small then they are conjugate by \cite[Theorem 2]{Prill-quotients}. The $\C^*$-ruling of $S_0$ does not extend to a ruling of $S'$, so by \cite[4.5]{Palka-classification1} its boundary is a rational fork with admissible maximal twigs and its minimal normal completion is unique up to isomorphism. (For the description of the boundary one could also refer to a more general result \cite[2.5.2.14]{Miyan-OpenSurf}.) \end{proof}

\subsection{Affine-ruled planes}\label{ssec:affine-ruled S'}

By \ref{prop4:non affine-ruled S'} we may assume that $S'$ is affine-ruled. This gives an affine ruling of $S_0$. We assume that $(\ov S,D+\E,p)$ is a minimal completion of the latter. This weakens our initial snc-minimality assumption on $D$, i.e. $D$ is now $p$-minimal, but the unique section contained in $D$ may be a non-branching $(-1)$-curve. The base of $p$ is rational, because it is isomorphic to a section contained in $D+\E$.

\begin{figure}[h]\centering\includegraphics[scale=0.4]{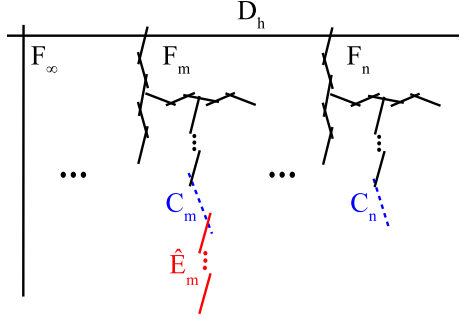}\caption{Affine-ruled $S'$}  \label{fig:affineruled}\end{figure}

\blem\label{lem4:affine-ruled S'} If $S'$ is affine-ruled then there exists exactly one fiber of $p$ contained in $D$ (see Fig.~\ref{fig:affineruled}). Each other singular fiber has a unique $(-1)$-curve, which is an $S_0$-component. The singularities of $S'$ are cyclic. \elem

\begin{proof} We have $\Sigma_{S_0}=\nu-1$  and $\nu\leq 1$ by \ref{prop:homology}, so $\Sigma_{S_0}=0$ and there is exactly one fiber $F_\8$ contained in $D$. The fiber is smooth by the $p$-minimality of $D$. Each singular fiber $F$ of $p$ contains exactly one $(-1)$-curve. Indeed, if $D_0\subseteq D$ is a vertical $(-1)$-curve then by the $p$-minimality of $D$ it meets $D_h$ and two $D$-components, so $\mu(D_0)>1$. The latter is impossible, as $D_h\cdot F=1$. The $(-1)$-curve, say $C$, has $\mu(C)>1$ and it is the unique $S_0$-component of $F$. There are exactly two components of multiplicity one in $F$, they are tips of $F$ and $D_h$ intersects one of them. Thus the connected component of $\un F-C$ not contained in $D$ is a chain, so $S'$ has only cyclic singularities. \end{proof}

\bsrem Note that in the lemma above it was pre-assumed (as in the whole paper, cf. section \ref{ssec:Qhp notation}) that $S'$ is logarithmic, but there is not need for this in fact. Namely, in any case $\E$ is vertical, so it is a rational forest. Then $D$ is a rational tree and $\ov S$ and the base of $p$ are rational by \cite[3.4(i)]{Palka-classification1}. The remaining part of the argument goes through. \esrem

\bcon\label{con4:affine-ruled S'} Let $\F_1=\PP(\cal O_{\PP^1}\oplus \cal O_{\PP^1}(-1))$ be the first Hirzebruch surface with a (unique) projection $\wt p\:\F_1\to \PP^1$. Denote the section coming from the inclusion of the first summand by $D_h'$, then $D_h'^2=-1$. Choose $n+1$ distinct points $x_\8,x_1,\ldots,x_n\in D_h'$ and let $F_\8$ be the fiber containing $x_\8$. For each $i=1,\ldots,n$ starting from a blow-up of $x_i$ create a fiber $F_i$ over $\wt p(x_i)$ containing a unique $(-1)$-curve $C_i$. Let $D_i$ be the connected component of $\un {F_i}-C_i$ intersecting $D_h$, the proper transform of $D_h'$. By renumbering we may assume there is $m\leq n$, such that $C_i$ is a tip of $F_i$ if and only if $i>m$. Assume also that $m\geq 1$ (for $m=0$ we would get a smooth surface). For $i\leq m$ put $\E_i= \un {F_i} -D_i -C_i$. Clearly, each $\E_i$ is a chain. Let $\ov S$ be the resulting surface and let $p\:\ov S\to \PP^1$ be the induced $\PP^1$-ruling. Put $D=F_\8+D_h+\sum_{i=1}^n D_i$, $S=\ov S-D$ and $\E=\sum_{i=1}^m\E_i$. We define $\epsilon\:S\to S'$ as the morphism contracting $\E_i$'s.\econ

\brem\label{rem3:H1 of affine-ruled S'} Let $p\:\ov S\to \PP^1$ be as in \ref{con4:affine-ruled S'} and for a fiber $F$ denote the greatest common divisor of multiplicities of all $S$-components of $F$ by $\mu_S(F)$. By \ref{prop:homology} $H_1(S',\Z)=H_1(S,\Z)$ and by \cite[4.19, 5.9]{Fujita} $$H_1(S,\Z)=\bigoplus_{i=1}^n\Z_{\mu_S(F_i)},$$ so $H_1(S',\Z)$ can be any finite abelian group. It is easy to see that $\mu_S(F_i)=\mu(C_i)/d(\E_i)$, where $d(\E_i)=d(0)=1$ if $i>m$. In particular, $S'$ is a $\Z$-homology plane if and only if $m=n$ and each $F_i$ is a chain. In fact then $\pi_1(S)$ vanishes, so $S'$ is contractible. \erem

\bthm\label{thm4:affine-ruled S'} The surface $S'$ constructed in \ref{con4:affine-ruled S'} is an affine-ruled singular $\Q$-homology plane. Conversely, each singular $\Q$-homology plane admitting an affine ruling can be obtained by construction \ref{con4:affine-ruled S'}. Its strongly balanced boundary is unique if it is branched and is unique up to reversion if it is a chain. The affine ruling of $S'$ is unique if and only if its strongly balanced boundary is not a chain.  \ethm

\begin{proof} By definition $\E_i$'s are admissible chains, so $S'$ is normal and has only cyclic singularities. We have $d(D)=-\prod_i d(D_i)$ (cf. \cite[2.1.1]{KR-C*onC3}), so $d(D)\neq 0$, hence $S'$ is a singular $\Q$-homology plane by \ref{lem5:when ruling gives sQhp}. The last part of the statement almost follows from \ref{lem4:affine-ruled S'}. It remains to note that by a flow (cf. \ref{ex2:flow}) we can change freely the self-intersection of the horizontal boundary component without changing the rest of $D$, so we can assume that the construction starts with a negative section on $\F_1$. (We could for instance start with $D_h'$ equal to the negative section on $\F_n$, so that the resulting boundary would be strongly balanced, cf. \ref{def2:normalization of boundary}). The uniqueness of a strongly balanced boundary follows from \ref{prop2:boundary uniqueness}.

We now consider the uniqueness of an affine ruling. Let $(V_i,D_i,p_i)$ be two minimal completions of two affine rulings of $S'$ (cf. \ref{ssec:rulings}). By \ref{lem4:affine-ruled S'} both $D_i$ contain a $0$-curve $F_{\8,i}$ as a tip. By flows with supports in $F_{\8,i}$ we may assume both $D_i$ are standard (cf. \ref{def2:standard boundary}).

Suppose $D_1$ is not a chain. Then $D_1$ and $D_2$ are isomorphic as weighted curves (cf. \ref{prop2:boundary uniqueness}). Let $T_i$ be the unique maximal twig of $D_i$ containing a $0$-curve. Either $T_i=F_{\8,i}=[0]$ or we can write $T_i=[0,0,a_1,\ldots,a_n]$ with $[a_1,\ldots,a_n]$ admissible. By \ref{prop2:boundary uniqueness} there is a flow $\Phi\:(V_1,D_1)\transf (V_2,D_2)$. Since $D_1$ is branched, $\Supp \Phi^\bullet\subseteq T_1$. Moreover, it follows from \ref{prop2:boundary uniqueness} and \ref{ex2:flow} that $\Supp \Phi^\bullet\subseteq F_{\8,i}$. For $i=1,2$ let $f_i$ be some fiber of $p_i$ different than $F_{\8,i}$. Since $\Phi^\bullet(f_1)$ is disjoint from $F_{\8,2}$, we get $\Phi^\bullet(f_1)\cdot f_2=0$, so $p_1$ and $p_2$ agree on $S'$.

Suppose now that $(V_1,D_1)$ is a standard completion of $S'$ with $D_1=[0,0,a_1,\ldots,a_n]$. We may assume that $[a_1,\ldots,a_n]$ is admissible and nonempty. Indeed, if it is empty then $S'\cong \C^2$ is smooth and if it is non-admissible then, by the Hodge index theorem we necessarily have $D_1=[0,0,0]$, which disagrees with \ref{prop:homology}(i). Let $(V_2,D_2)$ be another completion of $S'$ with $D_2$ being a reversion of $D$. The $0$-tip $T_i$ of each $D_i$ induces an affine ruling on $S'$. Let $(V,D)$ be a minimal normal pair dominating both $(V_i,D_i)$, such that both affine rulings extend to $\PP^1$-rulings of $V$. We argue that these affine rulings are different by proving that $\sigma_1^*T_1\cdot \sigma_2^*T_2\neq 0$, where $\sigma_i:(V,D)\to (V_i,D_i)$ are the dominations. Suppose $\sigma_1^*T_1\cdot \sigma_2^*T_2=0$. Let $H$ be an ample divisor on $V$ and let $(\lambda_1,\lambda_2)\neq (0,0)$ be such that $\wt T\cdot H=0$ for $\wt T=\lambda_1\sigma_1^*T_1+\lambda_2\sigma_2^*T_2$. We have $(\sigma_i^*T_i)^2=T_i^2=0$, so $$\wt T^2=2\lambda_1\lambda_2\sigma_1^*T_1\cdot \sigma_2^*T_2=0,$$ hence $\wt T\equiv 0$ by the Hodge index theorem. However, $D$ has a non-degenerate intersection matrix, because $d(D)=d(D_1)\neq 0$, so $\wt T$ is a zero divisor. Then $\sigma_1^*T_1=[0]$, otherwise $\sigma_1^*T_1$ and $\sigma_2^*T_2$ would contain a common $(-1)$-curve, which contradicts the minimality of $(V,D)$. It follows that $\sigma_1$ (and $\sigma_2$) are identities. This contradicts the fact that the reversion for nonempty $[a_1,\ldots,a_n]$ is a nontrivial transformation of the completion (even if $[a_1,\ldots,a_n]^t=[a_1,\ldots,a_n]$). \end{proof}

The following example shows that even if the strongly balanced boundary is unique, there might be infinitely many strongly balanced completions.

\bex\label{ex:one boundary many completions} Let $(V,D,\iota)$ be an snc-minimal completion ($\iota$ is the embedding, cf. \ref{ssec:divs and pairs}) of an affine-ruled singular $\Q$-homology plane $S'$ as above. Assume $D_h$ is branched and $D_h^2=-1$. The only change of $D$ which can be made by a flow is a change of the weight of $D_h$. If we now make an elementary transformation $(V,D)\mapsto (V_x,D_x)$ with a center $x\in F_\8\setminus D_h$ then $D$ becomes strongly balanced (cf. \ref{def2:normalization of boundary}). Denote the resulting completion by $(V_x,D_x,\iota_x)$ and let $F_{\8,x}$ be the new fiber at infinity. The isomorphism type of the weighted boundary $D_x$ does not depend on $x$, but the completions (triples) are clearly different for different $x$. Moreover, in general even the isomorphism type of the pair $(V_x,D_x)$ depends on $x$. To see this suppose $(V_x,D_x)\cong (V_y,D_y)$. As the isomorphism maps $F_{\8,x}$ to $F_{\8,y}$, we get an automorphism of $(V,D)$ mapping $x$ to $y$. Taking a minimal resolution $\ov S\to V$, contracting all singular fibers to smooth fibers without touching $D_h$ and then contracting $D_h$ we see that for $x\neq y$ this automorphism descends to a nontrivial automorphism of $\PP^2$ fixing points which are images of contracted $S_0$-components and of $D_h$. In general such an automorphism does not exist. \eex

\subsection{Moduli}\label{ssec:moduli}

Repeating the construction \ref{con4:affine-ruled S'} in a special case we will now obtain arbitrarily high-dimensional families of non-isomorphic singular $\Q$-homology planes with negative Kodaira dimension of the smooth locus and the same homeomorphism type. The following example gives a proof of the theorem \ref{thm:moduli}. For smooth $\Q$-homology planes a similar example was considered in \cite[4.16]{FZ-deformations}.

\begin{figure}[h]\centering\includegraphics[scale=0.4]{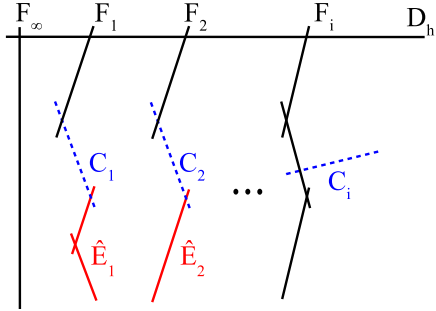}\caption{Singular fibers in example \ref{ex3:family of S' with fixed D}} \label{fig:moduli example}\end{figure}

\bex\label{ex3:family of S' with fixed D}  Put $m=2$ and $n=N+2$ for some $N>0$ and let $\ov S$, $D$, $\E$ etc. be created as in the construction above, so that $D_1=[3]$, $D_2=[2]$ and $D_i=[2,2,2]$ for $3\leq i\leq n$. Then $\E_1=[2,2]$ and $\E_2=[2]$ (see Fig.~\ref{fig:moduli example}).

Denoting the contraction of $\sum_{i=3}^n C_i$ by $\sigma:\ov S\to V$ we can factor the contraction $\ov S\to \F_1$ (which reverses the construction) as the composition $\ov S\xrightarrow{\sigma} V\xrightarrow{\sigma'}\F_1$. Put $y_i=\sigma(C_i)$ and $y=(y_3,\ldots,y_n)$. While $\sigma'^{-1}$ is determined uniquely by the choice of $(x_1,\ldots,x_n)$, $\sigma^{-1}$ and the resulting surface $\ov S$ (and hence $S'$) can depend on the choice of $y$. Let us write $\ov S_y$ and $S'_y$ to indicate this dependence. For $3\leq i\leq n$ let $D_i^0$ be the open subset of the middle component of $D_i$ remaining after subtracting two points belonging to other components of $D_i$. Put $$U=D_4^0\times\ldots\times D_n^0\cong \C^{N-1}.$$ The family $$\{S'_y\}_{y\in D_3^0\times U}\to D_3^0\times U$$ is $N$-dimensional. Since there exists a compactly supported auto-diffeomorphism of the pair $(\C^2,\C^*\times \{0\})$ mapping $(p,0)$ to $(q,0)$ for any $p,q\neq 0$, the choice of $y\in D_3^0\times U$ is unique up to a diffeomorphism fixing irreducible components of $\sigma_*(D+\E+C_1+C_2)$. Thus all $S'_y$ are homeomorphic.

Let $\pi:\mathfrak X\to U$ be the subfamily over $\{y_3^0\}\times U$. We will show that the fibers of $\pi$ are non-isomorphic. Suppose that $S_y'\cong S_z'$ for $y,z\in \{y_3^0\}\times U$. The isomorphism extends to snc-minimal resolutions. By \ref{prop2:boundary uniqueness} there is a flow $\Phi^\bullet:\ov S_y\map\ov S_z$, which is an isomorphism outside $F_\8$. Clearly, $\Phi^\bullet$ fixes $D_h\setminus\{x_\8\}$, $F_1$ and $F_2$, hence restricts to an identity on $D_h\setminus\{x_\8\}$ and respects fibers. Since $C_i$ are unique $(-1)$-curves of the fibers, they are fixed by $\Phi^\bullet$. It follows that $\Phi^\bullet_{|\ov S-F_\8-D_h}$ descends to an automorphism $\Phi_V$ of $V-F_\8-D_h$ fixing the fibers, such that $\Phi_V(y_i)=z_i$. Moreover, $\Phi_V$ descends to an automorphism $\Phi_{\F_1}$ of $\F_1-F_\8-D_h'$ fixing fibers. If $(x,y)$ are coordinates on $\F_1-F_\8-D_h'\cong \C^2$, such that $x$ is a fiber coordinate then $$\Phi_{\F_1}(x,y)=(x,\lambda y+P(x))$$ for some $P\in \C[x]$ and $\lambda\in\C$. Introducing successive affine maps for the blowups one can check that in some coordinates $\Phi_V$ acts on $D_i^0$ as $t\to\lambda^{\mu(C_i)} t$. Now the requirement $y_3=y_3^0$ fixes $\lambda^2=1$, so since $\mu(C_i)=2$ for each $3\leq i\leq n$, we get that $y=z$.\eex

\bsrem Note that by \cite[4.19, 5.9]{Fujita} for $S'$ as above $\pi_1(S')$ is the $N$-fold free product of $\Z_2$. It follows from \ref{rem3:H1 of affine-ruled S'} that given a weighted boundary there exist only finitely many affine-ruled singular $\Z$-homology planes with this boundary. That is why in the above example we have used branched fibers $F_i$ for $3\leq i\leq n$, so that the resulting surfaces are $\Q$-, but not $\Z$-homology planes. \esrem

\section{$\C^*$-ruled $\Q$-homology planes}\label{sec:classification of k(S') negative}

By \cite[1.1(2), 1.2]{Palka-classification1} and section \ref{ssec:affine-ruled S'} to accomplish the classification of singular $\Q$-homology planes with smooth locus of non-general type one needs to classify $\Q$-homology planes which are $\C^*$-ruled. Therefore in this section we assume that $S'$ is $\C^*$-ruled (and logarithmic, cf. \ref{ssec:Qhp notation}). The first homology group of $S'$ and some necessary conditions for singular fibers of such rulings have been analyzed in \cite[2.9 - 2.10]{MiSu-hPlanes}. As before, we concentrate on completions rather than the affine part itself, as this gives more information and allows to give a general method of construction. It also allows to compute the number of different $\C^*$-rulings, and in consequence the number of affine lines on $S'$.

\subsection{Properties of $\C^*$-rulings}\label{ssec:prop of C*-rulings}

We can lift the $\C^*$-ruling of $S'$ to a $\C^*$-ruling of the resolution and extend it to a $\PP^1$-ruling $p\:\ov S\to \PP^1$ of a smooth completion. Assume that $D+\E$ is $p$-minimal. By \ref{prop:homology}(v) $\Sigma_{S_0}=h+\nu-2$ and $\nu\leq 1$, so $(h,\nu,\Sigma_{S_0})= (1,1,0), (2,1,1)$ or $(2,0,0)$. The original $\C^*$-ruling of $S'$ is twisted with the base $\C^1$ in the first case, untwisted with the base $\C^1$ in the second case and untwisted with the base $\PP^1$ in the third case.

\blem\label{lem5:F_0 basic properties} Let $F_1,\ldots,F_n$ be all the columnar fibers of $p:\ov S\to \PP^1$ (cf. \ref{def2:(un)twisted and columnar fiber}). Let $F_\8$ be the fiber contained in $D$ if $\nu=1$. There is exactly one more singular fiber $F_0$, the fiber contains $\E$. We have also: \benum[(i)]

\item if $(h,\nu)=(1,1)$ then $F_\8=[2,1,2]$, $\sigma(F_0)=1$ and $F_0$ and $F_\8$ contain branching points of $p_{|D_h}$,

\item if $(h,\nu)=(2,1)$ then $F_\8$ is smooth and $\sigma(F_0)=2$,

\item if $(h,\nu)=(2,0)$ then $\sigma(F_0)=1$ and $F_0$ contains a $D$-rivet,

\item if $h=2$ then the components of $D_h$ are disjoint. \eenum \elem

\begin{proof} Suppose $(h,\nu)=(1,1)$. Then $\Sigma_{S_0}=0$, so by \cite[7.6]{Fujita} every singular fiber different than $F_\8$ is either columnar or contains a branching point of $p_{|D_h}$. Now $D_h$ is rational and $p_{|D_h}$ has two branching points, one of them contained in $F_\8$, as $D$ is a tree. It follows that $F_0$ is unique. The $p$-minimality of $D$ implies that $F_\8=[2,1,2]$. Assume now that $h=2$. We have $\Sigma_{S_0}=\nu\in \{0,1\}$ and the $p$-minimality of $D$ gives (ii), (iii) and the uniqueness of $F_0$. Suppose the components of $D_h$ have a common point. $D$ is a tree, so in this case $\nu=0$, which gives $\sigma(F_0)=1$. As $D$ is a simple normal crossing divisor, the common point belongs to the unique $S_0$-component of $F_0$, which has therefore multiplicity one. The connectedness of $D$ implies that $F_0$ contains no $D$-components. However, then $F_0$ has a unique $(-1)$-curve and it has multiplicity one, which is impossible by \ref{lem2:fiber with one exc comp}. \end{proof}

The above lemma is essentially the lemma \cite[2.10]{MiSu-hPlanes}. While the conditions stated above are necessary, they are not sufficient. In the following examples the $\C^*$-ruling satisfies 2.10 loc. cit. but the $\C^*$-ruled surface one obtains is not a $\Q$-homology plane.

\bex\label{ex: C*-ruled with d(D)=0:I} Let $\F_n$, $n\geq 0$, be the $n$-th Hirzebruch surface and let $D_0$, $D_\8$ be sections with $D_0^2=n$ and $D_\8^2=-n$. Let $F_\8$ be a fiber, put $D=D_0+D_\8+F_\8$. Pick a point not belonging to $D$ and make a connected sequence of blow-ups over it. Let $C_0$ be the unique $(-1)$-curve in the inverse image of the point and let $F_0$ and $C_1$ be the reduced total and the proper transform of the fiber. Denote the resulting surface by $\ov S$, put $S=\ov S-D$, $\E=F_0-C_0-C_1$ and let $S\to S'$ be the morphism contracting $\E$. In particular $\E$ can be any admissible chain and then $S'$ has a unique cyclic singular point. Note that $S'$ is not a $\Q$-homology plane because $d(D)=0$ (cf. \ref{lem5:when ruling gives sQhp}(iv)). \eex

\bex\label{ex: C*-ruled with d(D)=0:II} Start with the pair $(\F_1,D_0+D_\8)$, where $\F_1$ is the first Hirzebruch surface, $D_0$ and $D_\8$ are sections with $D_0^2=1$ and $D_\8^2=-1$. Pick two points on $D_0$ and blow up over it to create two singular fibers $F_1=[2,1,2]$, $F_2=[2,1,2]$. Denote their $(-1)$-curves by $C_1$, $C_2$. These $(-1)$-curves separate two chains  $T_0=[2,1,2]$ and $T_\8=[2,1,2]$, where the middle $(-1)$-curves are $D_0$ and $D_\8$ respectively. We have $d(T_0)=d(T_\8)=0$. Now pick a point on some $C_i$, say $C_1$, which does not belong to $T_0+T_\8$ and make a connected sequence of blow-ups over it. Let $C_0$ be the unique $(-1)$-curve in the inverse image of the point and let $F_0$ be the total reduced transform of the fiber. Denote the resulting complete surface by $\ov S$. If $C_0$ is not a tip of $F_0$ then denote the connected component of $F_0-C_0$ not meeting $D_0+D_\8$ by $\E$. Let $D$ be the reduced divisor with support $T_0\cup T_\8\cup (F_0-C_0-\E)$. Put $S=\ov S-D$, $\E=F_0-C_0-C_1$ and let $S\to S'$ be the morphism contracting $\E$ (which is necessarily an admissible chain). Once again $S'$ is not a $\Q$-homology plane because $d(D)=0$. \eex

Theoretically, if $X$ is a normal surface and $p':X\to B$ is a $\C^*$-ruling then taking a completion of $X$ and an  extension of $p'$ to a $\PP^1$-ruling we are able, using \ref{lem5:when ruling gives sQhp}, to recognize when $X$ is a $\Q$-homology plane (note that in particular $B$ has to be rational). However, to give constructions we need to reformulate the condition $d(D)\neq 0$ in a way which is easier to verify by looking at the geometry of singular fibers. Recall that for a family of subsets $(A_i)_{i\in I}$ of a topological space $Y$ a subset $X\subseteq Y$ \emph{separates} the subsets $(A_i)_{i\in I}$ (inside $Y$) if and only if each $A_i$ is contained in a closure of some connected component of $Y\setminus X$ and none of these closures contains more than one $A_i$. Recall also that by convention a twig of a fixed divisor is ordered so that its tip is the first component.

\blem\label{lem5:when d(D)=0} Let $(\ov S,T,p)$ be a triple satisfying conditions \ref{lem5:when ruling gives sQhp}(i)-(iii). Assume additionally that $T$ is $p$-minimal and $f\cdot T=2$ for a general fiber $f$ of $p$. In case $(h,\nu)=(2,0)$ let $D_0$, $F_0$, $B$, $\wt D_0$ be respectively some horizontal component of $D$, a unique fiber containing a $D$-rivet, a unique component of $D$ separating $D_0$, $D_h-D_0$ and $\E$ inside $D\cup F_0$ and a connected component of $D-B$ containing $D_0$. Then $d(D)\neq 0$ if and only if the following conditions hold:\benum[(i)]

\item the base of the fibration is $\PP^1$ or $\C^1$ (i.e. $\nu\leq 1$),

\item if $(h,\nu)=(2,1)$ then both $\ov S-T$-components of the fiber with $\sigma=2$ intersect $D$,

\item if $(h,\nu)=(2,0)$ then $d(\wt D_0)\neq 0$.

\eenum\elem

Note that the advantage of the condition (iii) over $d(D)\neq 0$ is that $\wt D_0$ is simpler than $D$, it contains at most one branching component.

\begin{proof} Clearly, if $d(D)\neq 0$ then $S'$ is a $\Q$-homology plane by \ref{lem5:when ruling gives sQhp}, which implies (i) and (ii) ($D$ meets each curve not contained in $D+\E$ because $S'$ is affine). Suppose now that the conditions (i) and (ii) are satisfied. We show that $d(D)\neq 0$ is equivalent to (iii) (which is an empty condition if $(h,\nu)\neq (2,0)$). Note that $d(D)\neq 0$ is equivalent to $d(T)\neq 0$, as $T-D$ is negative definite.

Consider the case $h=1$. We have $\Sigma_{\ov S-T}=\nu-1$, hence $\nu=1$ and $\Sigma=0$. The horizontal component $D_h$ meets the unique fiber $F_\8$ contained in $T$ in one point, because $T$ is a forest. Let $T_\8$ be the component meeting $D_h$. We have $d(F_\8)=0$, so by \cite[2.1.1(i)]{KR-C*onC3} $$d(D)=d(F_\8)d(D-F_\8)-d(F_\8-T_\8)d(D-F_\8-D_h)$$ and we obtain $$d(D)=-d(F_\8-T_\8)d(D-F_\8-D_h).$$ Since $F_\8-T_\8$ and $D-F_\8-D_h$ are vertical and do not contain whole fibers, they are negative definite, hence $d(D)<0$.

We may now assume $h=2$. Then $\Sigma=\nu\in \{0,1\}$. Put $\E=T-D$. In case $\nu=1$ let $F_\8$ be the unique fiber contained in $D$, and let $F_0$ be the unique singular fiber with $\sigma(F_0)=2$. In case $\nu=0$ let $F_0$ be the unique fiber containing a $D$-rivet. All other singular fibers are columnar by \cite[7.6]{Fujita}, so they contain no components of $\E$. We need to prepare some tools to proceed. Recall that the Neron-Severi group of $\ov S-T$ is defined as the quotient of $NS(\ov S)$ by the subgroup generated by components of $T$. We put $\rho(\ov S-T)=\dim NS(\ov S-T)\otimes \Q$.

Let $(X,R)$ be a smooth pair, $X$ rational. Suppose $R=R_1+R_2$, where $R_1$ and $R_2$ meet in unique components $C_1\subseteq R_1$, $C_2\subseteq R_2$ respectively. If at least one of $R_i$, $i=1,2$ is negative definite then we call $R-C_1$ \emph {a swap} of $R-C_2$ and vice versa. Similarly $(X,R-C_i)$, $i=1,2$ are by definition swaps of each other and so are $X-(R-C_i)$ for $i=1,2$. The basic property of this operation we will need is that $\rho(X-(R-C_1))=\rho(X-(R-C_2))$. To see this it is enough to show that $C_1, C_2$ to not belong to the subspace $V$ of $NS(X)\otimes \Q$ generated by components of $R_1-C_1+R_2-C_2$. By symmetry we can assume $R_2$ is negative definite. Suppose $C_1\in V$, write $C_1\equiv U_1+U_2$, where $U_i$ is in the subspace generated by components of $R_i-C_i$. Then $0=C_1\cdot U_2=U_1\cdot U_2+U_2^2=U_2^2$, hence $U_2\equiv 0$ by negative definiteness of $R_2$. Then $0<C_1\cdot C_2=U_1\cdot C_2=0$, a contradiction. Suppose $C_2\in V$ and write $C_2\equiv U_1+U_2$ as above. Then $(C_2-U_2)^2=(C_2-U_2)\cdot U_1=0$, so  $C_2\equiv U_2$ by negative definiteness of $R_2$. Then $0<C_1\cdot C_2=C_1\cdot U_2=0$, a contradiction. Thus swapping preserves $\rho$. The definition is of general use, but below we use only a special kind of swapping, namely (keeping the assumption that $R_2$ is negative definite) when $C_2$ is a $(-1)$-curve and it is absorbed into the boundary, i.e. we do the swap one way, changing $(X,R-C_2)$ to $(X,R-C_1)$.

Now start with $(\ov S,T)$ and perform interchangeably contractions of $(-1)$-curves in $F_0$ (and its images) which are non-branching components of the boundary and swaps absorbing vertical $(-1)$-curves in $F_0$ (and its images) into the boundary. Denote the resulting smooth pair by $(X,T')$. By the properties of swaps and blowing-ups the rank of the Neron-Severi group of the open part and the difference between $b_2$ of the complete surface and the number of components in the boundary remains constant. We see also that $T'$ is a rational forest. Now the crucial remark is that $d(T)=0$ if and only if $d(T')=0$. To see this we may assume that $(X,T')$ is simply a swap of $(\ov S,T)$ as above. Since the number of components of $T$ equals $b_2(\ov S)$, $d(T)\neq 0$ if and only if $\rho(\ov S-T)=0$, which is equivalent to $\rho(X-T')=0$ and then to $d(T')\neq 0$.

Consider the case $\Sigma=\nu=0$. At some point the process of swapping and contracting makes $B$ into a $0$-curve or a $(-1)$-curve. It is easy to see that the divisor $\wt D_0+\wt D_\8$ is not affected by the process, so we have $d(D)\neq 0$ if and only if $d(\wt D_0)\cdot d(\wt D_\8)\neq 0$. All singular fibers of the induced $\PP^1$-ruling at this stage are columnar, so they can be written as $R_{i,0}+C_i+R_{i,\8}$, where $i=1,\ldots,n'$ enumerates these fibers, $C_i^2=-1$, $R_{i,0}$ and $R_{i,\8}$ are chains whose last components meet $D_0$ and $D_\8$ respectively. For $j=0,\8$ put $\wt e_j=\wt e(\wt D_j)$ (cf. \ref{ssec:divs and pairs}). Then $\wt e_j=\sum_{i} \wt e(R_{i,j})$. We have $d(\wt D_j)=(-D_j^2-\wt e_j)\cdot \prod_{i}d(R_{i,j})$. By the properties of columnar fibers $$d(\wt D_0)+d(\wt D_\8)= -(D_0^2+D_\8^2+n')\cdot \prod_{i}d(R_{i,0}).$$ Moreover, when contracting singular fibers to smooth ones $D_0+D_\8$ is touched $n'$ times and its image consists of two disjoint sections on a Hirzebruch surface. It follows that $D_0^2+D_\8^2+n'=0$ and hence $d(\wt D_\8)+d(\wt D_0)=0$. Thus $d(D)\neq 0$ if and only if $d(\wt D_0)\neq 0$.

Consider the case $\Sigma=\nu=1$. We first show that $T'$ has at most one horizontal component. Suppose it has two. Then $\sigma(\wt F_0)=\sigma(F_0)=2$, so $\wt F_0$ contains a $(-1)$-curve, say $C_1$. Since $T'$ is $p$-minimal, $C_1\not\subseteq T$. Since we assumed that every $\ov S-T$ -component meets $D$, by the properties of swaps every $X-T'$ -component meets $T'$. By the definition of $X$ absorbing the $(-1)$-curve by a swap into the boundary is impossible. In particular if $\wt F_0$ has no more $(-1)$-curves then $C_1$ is not a tip of $\wt F_0$, so $\wt F_0$ is a chain. However, since $\sigma(\wt F_0)=2$, a swap absorbing $C_1$ into the boundary is possible, a contradiction. Thus $\wt F_0$ has two $(-1)$-curves, $C_1$ and $C_2$. One of them meets some horizontal component of $T'$, otherwise either $C_1$ or $C_1$ is a tip or $\wt F_0\cap T'$ has three connected components, in each case a swap absorbing one of $C_i$'s into the boundary would be possible. But then a similar argument shows that also the second $(-1)$-curve meets a horizontal component of $T'$. Thus $\wt F_0'$ is a chain with $C_1$, $C_2$ as tips and again a swap a before is possible, a contradiction. Thus $T'$ has at most one horizontal component. However, after the first swap where $\sigma$ of the image of $F_0$ drops the fiber has only one $(-1)$-curve, which has therefore multiplicity greater than one, hence no more swaps of this kind are possible. Thus $T'$ has a unique horizontal component $T'_h$. Then $$d(T')=d(F_\8)d(T'-F_\8)-d(T'-F_\8-D_\8)=-d(T'-F_\8-D_\8).$$ Now $T'-F_\8-D_\8$ is vertical and does not contain whole fibers, hence it is negative definite and we obtain $d(T')=d(T'-F_\8-D_\8)\neq 0$. \end{proof}

\bsrem By \ref{prop:homology} for any $\Q$-homology plane we have $H_i(S',\Z)=0$ for $i>1$ and $$|H_1(S',\Z)|^2= \frac{d(D)}{d(\E)},$$ hence $S'$ is a $\Z$-homology plane if and only if $d(D)=d(\E)$. For a $\C^*$-ruled $S'$ more explicit computations are done in \cite{MiSu-hPlanes}, so we will not repeat it here. For example, by 2.17 loc. cit. if a $\Z$-homology plane with $\ovk(S_0)\neq -\8$ is $\C^*$-ruled then $\ovk(S_0)=1$ and the ruling is untwisted with base $\PP^1$. The conditions for $S'$ having such a ruling to be contractible are given in 2.11. loc. cit (in particular $n=2$). \esrem

\subsection{The Kodaira dimension}\label{ssec:Kodaira dimensions}

In \cite[2.9 - 2.17]{MiSu-hPlanes} one can find formulas for the Kodaira dimension of the smooth locus, $\ovk(S_0)$ in terms of properties of singular fibers of the $\C^*$-ruling (the reader should note that in loc. cit. $\ovk(S')$ is by definition equal to $\ovk(S_0)$). Unfortunately the formulas 2.14(4), 2.15(2), 2.16(2) loc. cit. are incorrect. The corrections require splitting into cases depending on additional properties of singular fibers. We also compute the Kodaira dimension of $S'$. We keep the notation for singular fibers as in \ref{lem5:F_0 basic properties}. In case $\nu=0$ put $F_\8=0$. Let $J$ be the reduced divisor with the support equal to $D\cup F_0$. For $i=1,\ldots,n$ denote the $(-1)$-curve of the columnar fiber $F_i$ by $C_i$ and the multiplicity of $C_i$ by $\mu_i$. Put $J^+=J+C_1+\ldots+ C_n$.

\blem\label{lem5:minimal model of a Qhp} The divisor $J^+$ has simple normal crossings. Contract vertical $(-1)$-curves in $J^+$ and its images as long as the image is an snc-divisor. Let $\zeta\:(\ov S,J^+)\to (W,\zeta_*J^+)$ be the composition of these contractions. Then $\zeta_*F_i$ are smooth for $i=1,\ldots,n$ and:\benum[(i)]

\item if $h=1$ then $\zeta_*F_0=[2,1,2]$, $(\zeta_*D_h)^2=0$ and one can further contract $\zeta_*F_0$ and $F_\8$ to smooth fibers so that $W$ maps to $\F_1$ and $\zeta_*D_h$ maps to a smooth $2$-section of the $\PP^1$-ruling of $\F_1$ disjoint from the negative section,

\item if $h=2$ then $\zeta_*F_0$ is smooth, $W$ is a Hirzebruch surface and the components of $\zeta_*D_h$ are disjoint. Moreover, at least one of the components of $D_h$ has negative self-intersection and changing $\zeta$ if necessary one can assume that it is not affected by $\zeta$. \eenum\elem

\begin{proof} Suppose that $J^+$ does not have normal crossings at $x$. By \ref{lem5:F_0 basic properties} this can happen only if $h=2$. Moreover, $x\in D_h\cap F_0$, it is a branching point of $p_{|D_h}$ and two components of $F_0$ of multiplicity one meet at $x$. As $D$ has normal crossings, one of them is the unique $S_0$-component of $F_0$. By the $p$-minimality of $D$ it has to be a unique $(-1)$-curve of $F_0$ too, which is impossible by \ref{lem2:fiber with one exc comp}(i). This shows that $J^+$ is an snc-divisor. Since $F_i$ for $i=1,\ldots,n$ are columnar, $\zeta_*F_i$ are smooth.

Suppose $h=2$. Write $D_h=H+H'$. By \ref{lem5:F_0 basic properties} $H$ and $H'$ are disjoint. Since $H$ and $H'$ meet $F_0$ only in the components of multiplicity one, it follows from the definition of $\zeta$ that the images of $H'$ and $H$ intersect the same component of $\zeta_*F_0$. But this is possible only if $\zeta_*F_0$ is smooth. Since $\zeta_*J^+$ is snc, these images are disjoint. Say $H'^2\leq H^2$. Choosing the contracted $(-1)$-curves correctly we may assume that $H'$ is not affected by $\zeta$. Since $\zeta_*D_h$ consists of two disjoint sections on a Hirzebruch surface, we have $(\zeta_*D_h)^2=0$, so $D_h^2\leq 0$. Suppose $H^2=H'^2=0$. Then $\zeta$ does not affect $D_h$, so $n=0$ and $H$ and $H'$ intersect the same component $B$ of $F_\8$. If $\nu=1$ then $B$ is an $S_0$-component and the second $S_0$-component of $F_0$ does not intersect $D$, a contradiction with the affiness of $S'$. Thus $\nu=0$ and the condition \ref{lem5:when d(D)=0} is not satisfied (in other words $d(D)=0$), a contradiction.

Suppose $h=1$. By the definition of $\zeta$ the image of $D_h$ intersects the unique $(-1)$-curve of $\zeta_*F_0$. It follows that $\zeta_*F_0=[2,1,2]$. Now after the contraction of $F_0$ and $F_\8$ to smooth fibers the image of $W$ is a Hirzebruch surface $\F_N$, where $N\geq 0$, and the image $D_h'$ of $D_h$ is a smooth 2-section. Write $D_h'\equiv \alpha f+2H$ where $H$ is a section with $H^2=-N$ and $f$ is a fiber of the induced $\PP^1$-ruling of $\F_N$. We compute $$p_a(\alpha f+2H)=\alpha-N-1,$$ so since $D_h'$ is smooth, its arithmetic genus vanishes and then $\alpha=N+1$. Moreover, $D_h'\cdot H=\alpha-2N$, hence $D_h'\cdot H+N=1$. Now if $N=0$ then $\F_N=\PP^1\times \PP^1$ and an elementary transformation with center equal to the point of tangency of $D_h'$ and the image of $F_\8$ (which corresponds to a different choice of components to be contracted in $F_\8$) leads to $N=1$ and $D_h'\cdot H=0$. \end{proof}

\brem\label{rem5:blowing K+D} Let $(X,D)$ be a smooth pair and let $L$ be the exceptional divisor of a blow-up $\sigma:X'\to X$ of a point in $D$. Then $$K_{X'}+\sigma^{-1}D=\sigma^* (K_X+D)$$ if $\sigma$ is subdivisional for $D$ and $$K_{X'}+\sigma^{-1}D=\sigma^* (K_X+D)+L$$ if $\sigma$ is sprouting for $D$. \erem

Decompose $\zeta$ into a sequence of blow-downs $\zeta=\sigma_k\circ\ldots\circ \sigma_1$ and let $m\leq k$ be such that for $j>m$ the blow-up $\sigma_j$ is subdivisional for $(\sigma_j\circ\ldots\circ \sigma_1)_*J^+$ and $m$ is minimal such. Define $\eta\:\ov S\to \wt S$ and $\theta\:\wt S\to W$ as $\eta=\sigma_m\circ\ldots\circ \sigma_1$ and $\theta=\sigma_k\circ\ldots\circ \sigma_{m+1}$. Clearly, $\eta$ is an identity outside $F_0$. We denote a general fiber of a $\PP^1$-ruling by $f$.

\blem\label{lem5:K+eta_*J=} Let $\eta:\ov S\to \wt S$ and $\theta\:\wt S\to W$ be as above. Then $$K_{\wt S}+\eta_*J\equiv (n+\nu-1-\sum_{i=1}^n \frac{1}{\mu_i}) f+G+\theta^*\frac{1}{2}(U+U'),$$ where $G$ is a negative definite effective divisor with the support contained in $\Supp (F_\8+\sum_{i=1}^n F_i)$ and where $U$, $U'$ are the $(-2)$-tips of $\zeta_*F_0$ in case $p$ is twisted and are zero otherwise.\elem

\begin{proof} Let $V\subseteq W$ be defined as the sum of (four) $(-2)$-tips of $\un {F_\8}+\zeta_*\un {F_0}$ if $p$ is twisted and as zero otherwise. We check easily that $$K_{W}+D_h+\un {F_\8}+\zeta_*\un {F_0}\equiv (\nu-1)f+\frac{1}{2}V.$$ Indeed, if $p$ is untwisted this is just $K_{W}+D_h+2f\equiv 0$ on a Hirzebruch surface and if $p$ is twisted then it follows from the numerical equivalences $K_{W}+D_h+f\equiv 0$ and $\un {F_\8}+\zeta_*\un {F_0}-\frac{1}{2}V\equiv f$. By \ref{rem5:blowing K+D} we get $$K_{\wt S}+\eta_*J^+\equiv (n+\nu-1) f+\theta^*\frac{1}{2}V.$$ For every $i=1,\ldots,n$ the divisor $G_i=\frac{1}{\mu_i}F_i-C_i$ is effective and negative definite, as $C_i$ is not contained in its support. We get $$K_{\wt S}+\eta_*J\equiv (n+\nu-1) f+\sum_{i=1}^n(G_i-\frac{1}{\mu_i}F_i)+\theta^*\frac{1}{2}V,$$ so $$K_{\wt S}+\eta_*J\equiv (n+\nu-1-\frac{1}{\mu_i})f +\sum_{i=1}^nG_i+\theta^*\frac{1}{2}V$$ and we are done. \end{proof}

\brem\label{rem5:k and k_0 determine Kod dim} Since $K_{\ov S}+D+\E$ and $K_{\ov S}+D$ intersect trivially with a general fiber, we can write $K_{\ov S}+D+\E\equiv \kappa_0 f+G_0$ and $K_{\ov S}+D+\E\equiv \kappa f+G$, where $G_0$ and $G$ are some vertical effective and negative definite divisors and $\kappa_0,\kappa \in \Q$. It follows that $\ovk(S_0)$ and $\ovk(S)$ are determined by the signs of numbers $\kappa_0$ and $\kappa$ respectively. More explicitly, $\ovk(S_0)=-\8$, $0$, $1$ depending wether $\kappa_0<$, $=$ or $>0$ respectively. Analogous remarks hold for $\ovk(S)$ and $\kappa$. \erem

It turns out that $\kappa$ and $\kappa_0$ depend in a quite involved way on the structure of $F_0$. This dependence can be stated in terms of the properties of $\eta:\ov S\to\wt S$ defined above. Denote the $S_0$-components of $F_0$ by $C$, $\wt C$ (or just $C$ if there is only one) and their multiplicities by $\mu$, $\wt \mu$ respectively. Note that $\mu\geq 2$ if $\sigma(F_0)=1$, but if $\sigma(F_0)=2$ then it can happen that $\mu=1$ or $\wt \mu=1$.

\bthm\label{thm5:(k,k0) computations} Let $\lambda=n+\nu-1-\sum_{i=1}^n \frac{1}{\mu_i}$. The numbers $\kappa$ and $\kappa_0$ determining the Kodaira dimension of a $\C^*$-ruled singular $\Q$-homology plane $S'$ and of its smooth locus $S_0$ defined in \ref{rem5:k and k_0 determine Kod dim} are as follows:\benum[(A)]

\item Case $(h,\nu)=(1,1)$. Denote the component of $F_0$ intersecting the 2-section contained in $D$ by $B$.\benum[(i)]

   \item If $\eta=id$ and $F_0=[2,1,2]$ then $\kappa=\kappa_0=\lambda-\frac{1}{2}$.

   \item If $\eta=id$, $B$ is not a tip of $F_0$ and $C\cdot B>0$ then $(\kappa,\kappa_0)=(\lambda-\frac{1}{2},\lambda-\frac{1}{2\mu})$.

   \item If $\eta=id$, $C\cdot B=0$ and $F_0$ is a chain then $(\kappa,\kappa_0)=(\lambda-\frac{1}{2},\lambda)$.

   \item If $\eta=id$ and $B$ is a tip of $F_0$ then $(\kappa,\kappa_0)=(\lambda-\frac{1}{2},\lambda-\frac{1}{\mu})$.

  \item If $\eta\neq id$ then $\kappa=\kappa_0=\lambda$. \eenum

\item Case $(h,\nu)=(2,1)$. \benum[(i)]

   \item If $\eta=id$ and $C^2=\wt C^2=-1$ then $(\kappa,\kappa_0)=(\lambda-1,\lambda-\frac{1}{min(\mu,\wt \mu)})$.

   \item If $\eta=id$ and $C^2\neq -1$ or $\wt C^2\neq -1$ then $\kappa=\kappa_0=\lambda-\frac{1}{min(\mu,\wt \mu)}$.

   \item If $\eta\neq id$ then, assuming that $C$ is the $S_0$-component disjoint from $\E$, $\kappa=\kappa_0=\lambda-\frac{1}{\mu}$. \eenum

\item Case $(h,\nu)=(2,0)$. Then $\kappa=\kappa_0=\lambda$. \eenum \ethm

\begin{proof} (A) The unique $S_0$-component $C$ of $F_0$ is a $(-1)$-curve. Indeed, otherwise the $p$-minimality of $D$ implies that $B$ is the only $(-1)$-curve in $F_0$ and it intersects two other $D$-components of $F_0$, which gives $F_0=[2,1,2]\subseteq D$ with no place for $C$. It is now easy to check that the list of cases in (A) is complete. As $C^2=-1$, $\un {F_0}-C$ has at most two connected components. We see also that the only case where $\E$ is not connected is when $F_0$ contains no $D$-components, which is possible only if $C=B$ and $F_0=[2,1,2]$. Since $C$ is the unique $(-1)$-curve in $F_0$, $\zeta=\theta\circ\eta$ has at most one center on $\zeta_*F_0$, so by symmetry we can and will assume that it does not belong to $U'$ (cf. \ref{lem5:K+eta_*J=}). Suppose $\eta\neq id$. The center of $\eta$ belongs to a unique component of $\eta_*J$ and $D_h$ does not intersect components contracted by $\eta$. Then the mentioned component is a proper transform of a $D$-component, so $\eta_*(C+\E)=0$ by the connectedness of $\E$. If we now factor $\eta$ as $\eta=\sigma\circ\eta'$, where $\sigma$ is a sprouting blow-up for $\eta_*J$ then by \ref{lem5:K+eta_*J=} and \ref{rem5:blowing K+D} we get $$K+\sigma^{-1}\eta_*J\equiv \lambda f+G+\sigma^*\theta^*\frac{1}{2}(U+U')+Exc(\sigma),$$ where $Exc(\sigma)$ is the exceptional $(-1)$-curve contracted by $\sigma$ and $K$ is a canonical divisor on a respective surface. Since $\eta_*(C+\E)=0$, each component of $C+\E$ will appear with positive integer coefficient in $\eta'^*Exc(\sigma)$, which leads to $K_{\ov S}+\eta^{-1}\eta_*J\equiv \lambda f+G+G_0$, where $G_0$ is a vertical effective and negative definite divisor for which $G_0-\E-C$ is still effective. Since $\eta^{-1}\eta_*J=J=D+\E+C$, we get $\kappa=\kappa_0=\lambda$. We can now assume that $\eta=id$, so $$K_{\ov S}+D+\E+C\equiv \lambda f+G+\frac{1}{2}(U'+\theta^*U).$$ The latter can be written as $$K_{\ov S}+D\equiv (\lambda-\frac{1}{2}) f+G+\frac{1}{2}(U'+F_0+\theta^*U-2C-2\E).$$ All components of $F_0$ appear in $U'+F_0+\theta^*U$ with coefficients bigger than $1$, so $U'+F_0+\theta^*U-2C-2\E$ is effective and negative definite, as its support does not contain the $\E$-component which is a proper transform of $U$. This gives $\kappa=\lambda-\frac{1}{2}$. We now compute $\kappa_0$. If $F_0=[2,1,2]$ then $\theta^*U=U$ and $\E=U+U'$, so $K_{\ov S}+D\equiv (\lambda-\frac{1}{2})f+G$ and we get $\kappa_0=\lambda-\frac{1}{2}$. Suppose $B$ is a tip of $F_0$. Since $\mu(B)=2$, $F_0$ is a fork with two $(-2)$-tips as maximal twigs (cf. \ref{lem2:fiber with one exc comp}(ii)) and $\theta^*U=U$ ($U$ and $U'$ are components of $\E$). The divisor $G_0=\frac{1}{2}(U+U')+\frac{1}{\mu}F_0-C$ is vertical effective and its support does not contain $C$. Writing $$K_{\ov S}+D+\E\equiv (\lambda-\frac{1}{\mu})f+G+G_0$$ we infer that $\kappa_0=\lambda-\frac{1}{\mu}$, hence we obtain (iv). Consider the case (ii). Since $B$ is not a tip of $F_0$, $F_0$ is a chain. The assumption $B\cdot C>0$ implies that $B^2\neq -1$ and $\theta^*U=C+\E$. We obtain $$K_{\ov S}+D+\E\equiv(\lambda-\frac{1}{2\mu})f +G +\frac{1}{2}(U'+\E+\frac{1}{\mu}F_0-C)$$ and $U'+\E+\frac{1}{\mu}F_0-C$ is effective with support not containing $C$. This gives $\kappa_0=\lambda-\frac{1}{2\mu}$. We are left with the case (iii). As in (ii) $F_0$ is a chain and we have now $$K_{\ov S}+D+\E\equiv \lambda f+G+\frac{1}{2}(U'+\theta^*U-2C).$$ Since $B\cdot C=0$, $U'+\theta^*U-2C$ is effective and does not contain $B$, so $\kappa_0=\lambda$.

(B) Suppose $\eta\neq id$. Note that $\eta_*F_0$ contains a proper transform of one of $C$, $\wt C$, otherwise $F_0$ would contain a $D$-rivet. It follows that $\eta$ is a connected modification and its center lies on a birational transform of a $D$-component (the $S_0$-component contracted by $\eta$ has to intersect $D$). Thus $\eta_*F_0$ is a chain intersected by $D_h$ in two different tips and containing $C$. Since $D\cap \E=\emptyset$, we get $\eta_*(\wt C+\E)=0$. Writing $\eta=\sigma\circ\eta'$, where $\sigma$ is a sprouting blow-down, we see that $\eta'^*Exc(\sigma)$ is an effective negative definite divisor which does not contain $C$ in its support and for which $\eta'^*Exc(\sigma)-\wt C-\E$ is effective. By \ref{lem5:K+eta_*J=} we have $$K+\sigma^{-1}\eta_*D+C\equiv \lambda f+G+Exc(\sigma),$$ where $K$ is a canonical divisor on a respective surface. It follows from \ref{rem5:blowing K+D} and from arguments analogous to these from part (A) that $\kappa=\kappa_0=\lambda-\frac{1}{\mu}$. We can now assume that $\eta=id$. By \ref{lem5:K+eta_*J=} $$K_{\ov S}+D+C+\E+\wt C\equiv \lambda f+G,$$ which implies $\kappa_0=\lambda- \frac{1}{min(\mu,\wt \mu)}$. Writing $$K_{\ov S}+D\equiv (\lambda-\frac{1}{\alpha})f +G+\frac{1}{\alpha}(F_0-\alpha(C+\E+\wt C))$$ we see that $\kappa=\lambda-\frac{1}{\alpha}$, where $\alpha$ is the lowest multiplicity of a component of $C+\E+\wt C$ in $F_0$. Note that $C+\E+\wt C$ is a chain. Now if for instance $C^2\neq -1$ then $F_0$ is columnar and factoring $\theta$ into blow-downs we see that $\E$ is contracted before $C$, hence $\alpha=\mu\leq \wt \mu$. Suppose $C^2=\wt C^2=-1$ and let $\theta'$ be the composition of successive contractions of $(-1)$-curves in $F_0$ different than $C$. Now either $\theta'_*F_0=\theta'_*C=[0]$ or $\theta'_*F_0$ is columnar. Both possibilities imply that $C+\E$ contains a component of multiplicity one, hence $\alpha=1$.

(C) $C$ is a $(-1)$-curve. Indeed, $D\cap F_0$ contains at most one $(-1)$-curve and if it does then by the $p$-minimality of $D$ it meets both components of $D_h$ and has multiplicity one, so there is another $(-1)$-curve in $F_0$. We infer that $\un F_0-C$ has two connected components, one is $\E$ and the second contains a rivet. The existence of a rivet in $F_0$ implies that $\eta\neq id$, so $\eta_*(C+\E)=0$. Factoring out a sprouting blow-down from $\eta$ as above we get $$K+\sigma^{-1}\eta_*D \equiv \lambda f+G+Exc(\sigma).$$ The divisor $\eta'^*Exc(\sigma)-C-\E$ is effective and does not contain all components of $F_0$, so by \ref{rem5:blowing K+D} $\kappa=\kappa_0=\lambda$.  \end{proof}

\bsrem In case (B)(iii) it is not true in general that $\mu=min(\mu,\wt \mu)$.  \esrem

\subsection{Smooth locus of Kodaira dimension zero}\label{ssec:k(S_0)=0} As a corollary we obtain the following information in case $\ovk(S_0)=0$.

\bcor\label{cor5:when k_0=0} Let $S'$ be a $\C^*$-ruled singular $\Q$-homology plane and let $D$ be a $p$-minimal boundary for an extension $p$ of this ruling to a normal completion as above. Let $D$ be the $p$-minimal boundary and let $n$ be the number of columnar fibers. Then $\ovk(S_0)=0$ exactly in the following cases:\benum[(i)]

\item $n=0$ and $F_0$ is of type (A)(iii) or (A)(v),

\item $n=1$, $\mu=\mu_1=2$, $F_0$ contains no $D$-components and is of type (A)(i) or (A)(iv),

\item $p$ is untwisted with base $\C^1$, $n=1$, $\mu_1=2$, $min(\mu,\wt \mu)=2$ and some connected component of $F_0\cap D$ is a $(-2)$-curve,

\item $p$ is untwisted with base $\C^1$, $n=2$, $\mu_1=\mu_2=2$, and some $S_0$-component of $F_0$ meets $D_h$,

\item $p$ is untwisted with base $\PP^1$, $n=2$ and $\mu_1=\mu_2=2$. \eenum \ecor

\begin{proof} Note that $n-\sum_{i=1}^n\frac{1}{\mu_i}\geq \frac{n}{2}$, because $\mu_i\geq 2$ for each $i$. Suppose $p$ is twisted. Then $\mu\geq 2$, so by \ref{thm5:(k,k0) computations} $$\lambda\geq \kappa_0\geq \lambda-\frac{1}{2}\geq\frac{n-1}{2}.$$ If $n=0$ then $\lambda=0$, which gives $\kappa_0=0$ exactly in cases (A)(iii) and (A)(v). If $n=1$ then $\kappa_0=\lambda-\frac{1}{2}=0$, which is possible in case (A)(i) if $\mu_1=2$ and in case (A)(iv) if $\mu=\mu_1=2$. In both cases $D_h$ meets the $S_0$-component, so $F_0$ contains no $D$-components. If $p$ is untwisted with base $\PP^1$ then $$n-1\geq\lambda=\kappa_0\geq \frac{n}{2}-1,$$ so $n=2$ ($\lambda=-\frac{1}{\mu_1}<0$ for $n=1$) and $\kappa_0=1-\frac{1}{\mu_1}-\frac{1}{\mu_2}$, which vanishes only if $\mu_1=\mu_2=2$. Assume now that $p$ is untwisted with base $\C^1$. We have $$n>\kappa_0\geq \lambda-1\geq \frac{n}{2}-1,$$ so $n\in\{1,2\}$. There are no $(-1)$-curves in $D\cap F_0$ by the $p$-minimality of $D$, so at least one $S_0$-component, say $C$, is a $(-1)$-curve. We can also assume that $C$ is contracted by $\eta$ in case $\eta\neq id$ and that $\mu\geq \wt \mu$ in case $\eta=id$. Then $\kappa_0=\lambda-\frac{1}{\wt \mu}$. The composition $\xi$ of successive contractions of all $(-1)$-curves in $\un F_0-\wt C$ and its images is a connected modification. Suppose $n=2$. The inequalities above give $\lambda=1$, so $\mu_1=\mu_2=2$ and $\wt \mu=1$. Then $\xi_*F_0=[0]$ and since $\xi$ is a connected modification, $\wt C$ is a tip of $F_0$. It follows that some of $C$, $\wt C$ intersects $D_h$, otherwise $F_0-\wt C-C-\E$ is connected and intersects both sections from $D_h$, hence $F_0$ would contain a rivet. This gives (iv). Suppose $n=1$. Then $\mu_1=\wt \mu=2$. Note that by the choice of $C$ further contractions of $F_0$ to a smooth fiber are subdivisional for $\xi_*D\cup \xi_*F_0$, so $\xi_*F_0=[2,1,2]$ with the birational transform of $\wt C$ in the middle and the image of $D_h$ intersects both $(-2)$-tips of $\xi_*F_0$. Since $\xi$ is a connected modification, it does not touch one of these tips, so one of the connected components of $D\cap F_0$ is a $(-2)$-curve. Now if $\mu=1$ then $\mu<\wt \mu$, so by our assumption $\eta\neq id$. But then $\mu>1$, because $C^2=-1$ and $C$ intersects $\E$ and $D$. This contradiction ends the proof of (iii).\end{proof}

\subsection{Constructions}\label{ssec:construction}

Lemmas \ref{lem5:minimal model of a Qhp} and \ref{lem5:when ruling gives sQhp} give a practical method of reconstructing all $\C^*$-ruled $\Q$-homology planes. We summarize it in the following discussion. We denote irreducible curves and their proper transforms by the same letters.

\bcon \begin{case} \emph{A twisted ruling}. Let $D_h, x_0, x_\8$ be a smooth conic on $\PP^2$ and a pair of distinct points on it. Let $L_0$, $L_\8$ be tangents to $D_h$ at $x_0$, $x_\8$ respectively and let $L_i$ for $i=1,\ldots,n$, $n\geq 0$ be different lines (different than $L_0$, $L_\8$) through $L_0\cap L_\8$. Blow up at $L_0\cap L_\8$ once and let $p\:\F_1\to \PP^1$ be the $\PP^1$-ruling of the resulting Hirzebruch surface. Over each of $p(L_0)$, $p(L_\8)$ blow up on $D_h$ twice creating singular fibers $\wt F_0=[2,1,2]$ and $F_\8=[2,1,2]$. For each $i=1,\ldots, n$ by a connected sequence of blow-ups subdivisional for $L_i+D_h$ create a column fiber $F_i$ over $p(L_i)$ and denote its unique $(-1)$-curve by $C_i$. By some connected sequence of blow-ups with a center on $\wt F_0$ create a singular fiber $F_0$ and denote the newly created $(-1)$-curve by $C$ (if the sequence is empty define $C$ as the $(-1)$-curve of $\wt F_0$). Denote the resulting surface by $\ov S$, put $$T=D_h+\un {F_\8}+(\un {F_1}-C_1)+\ldots+(\un {F_n}-C_n)+\un {F_0}-C$$ and construct $S'$ as in \ref{lem5:when ruling gives sQhp}. $S'$ is a $\Q$-homology plane (singular if only $T$ is not connected), because conditions \ref{lem5:when ruling gives sQhp}(i)-(iii) are satisfied by construction and (iv) by \ref{lem5:when d(D)=0}. To see that each $S'$ admitting a twisted $\C^*$-ruling can be obtained in this way note that by the $p$-minimality of $D$ even if $F_0$ contains two $(-1)$-curves $C$ and $B\subseteq D$ then $B$ is not a tip of $F_0$ and $\zeta$ does not touch it, so in each case the modification $F_0\to \zeta_*F_0$ induced by $\zeta$ is connected and we are done by \ref{lem5:minimal model of a Qhp}. \end{case}

\begin{case}\emph{An untwisted ruling with base $\C^1$}. Let $x_0,x_1\ldots x_n,x_\8,y\in \PP^2$, $n\geq 0$ be distinct points, such that all besides $y$ lie on a common line $D_1$. Let $L_i$ be a line through $x_i$ and $y$. Blow up $y$ once and let $D_2$ be the negative section of the $\PP^1$-ruling of the resulting Hirzebruch surface $p\:\F_1\to \PP^1$. For each $i=0,1,\ldots,n$ by a connected sequence of blow-ups (which can be empty if $i=0$) with the first center $x_i$ and subdivisional for $D_1+L_i$ create a column fiber $F_i$ ($\wt F_0$ if $i=0$) over $p(x_i)$ and denote its unique $(-1)$-curve by $C_i$ if $i\neq 0$ and by $\wt C$ if $i=0$ (put $\wt C=L_0$ if the sequence over $p(x_0)$ is empty). Choose a point $z\in F_0$ which lies on $D_1+\un {\wt F_0}-\wt C$ and by a nonempty connected sequence of blow-ups with the first center $z$ create some singular fiber $F_0$ over $p(x_0)$, let $C$ be the new $(-1)$-curve. Denote the resulting surface by $\ov S$, put $$T=D_1+D_2+L_\8+(\un {F_1}-C_1)+\ldots+(\un {F_n}-C_n)+\un {F_0}-C-\wt C$$ and construct $S'$ as in \ref{lem5:when ruling gives sQhp}. The surface $S'$ is a $\Q$-homology plane by \ref{lem5:when d(D)=0}, as \ref{lem5:when d(D)=0}(ii) is satisfied by the choice of $z$. To see that all $S'$ admitting an untwisted $\C^*$-ruling with base $\C^1$ can be obtained in this way note that changing the completion of $S'$ by a flow if necessary we can assume that one of the components of $D_h$ is a $(-1)$-curve. Note also that, $D\cap F_0$ contains no $(-1)$-curves and, as it was shown in the proof of \ref{thm5:(k,k0) computations}, $\eta$ contracts at most one of $C$, $\wt C$. Then we are done by \ref{lem5:minimal model of a Qhp}. \end{case}

\begin{case}\emph{An untwisted ruling with base $\PP^1$}. Let $D_2$ be the negative section of the $\PP^1$-ruling of a Hirzebruch surface $p\:\F_N\to \PP^1$, $N>0$. Let $x_0,x_1,\ldots,x_n$, $n\geq 0$ be points on some section $D_1$ of $p$ disjoint from $D_2$. For each $i=0,1,\ldots,n$ by a connected sequence of blow-ups (which can be empty if $i=0$) with the first center $x_i$ and subdivisional for $D_1+p^{-1}(p(x_i))$ create a column fiber $F_i$ ($\wt F_0$ if $i=0$) over $p(x_i)$ and denote its unique $(-1)$-curve by $C_i$ if $i\neq 0$ and by $B$ if $i=0$ (put $B=p^{-1}(p(x_0))$ if the sequence over $p(x_0)$ is empty). Assume that the intersection matrix of at least one of two connected components of $$D_1+D_2+(\un {F_1}-C_1)+\ldots +(\un {F_n}-C_n)+(\un {\wt F_0}-B)$$ is non-degenerate. By a connected sequence of blow-ups starting from a sprouting blow-up for $D_1+\wt F_0$ with center on $B$ create some singular fiber $F_0$ over $p(x_0)$, let $C$ be the new $(-1)$-curve. Denote the resulting surface by $\ov S$, put $$T=D_1+D_2+(\un {F_1}-C_1)+\ldots+(\un {F_n}-C_n)+(\un {F_0}-C)$$ and construct $S'$ as in \ref{lem5:when ruling gives sQhp}. Note that $D$ is connected, because the modification $F_0+D_1\to \wt F_0+D_1$ is not subdivisional, so $S'$ is a $\Q$-homology plane by \ref{lem5:when d(D)=0}. By \ref{lem5:minimal model of a Qhp} and \ref{lem5:when d(D)=0} each $S'$ with an untwisted $\C^*$-ruling having a base $\PP^1$ can be obtained in this way. \end{case}  \econ

\section{Corollaries}

\subsection{Completions and singularities}\label{ssec:completion and singularities}

Recall that $\Q$-homology planes with non-quotient singularities have unique snc-minimal completions (and hence also the balanced ones) and unique singular points (cf. \cite{Palka-classification1}). The completions and singularities in case $\ovk(S_0)=-\8$ are described in section \ref{sec:k(S0) negative}. In case $\ovk(S_0)=2$ it is known (see loc. cit.) that the singular point is unique and of quotient type. Moreover, the snc-minimal boundary cannot contain non-branching $b$-curves with $b\geq 0$ as these induce $\C^1$ or $\C^*$-rulings of $S_0$, hence also the snc-minimal completion is unique. The theorem \ref{thm5:singularities}  summarizes information in the remaining cases.

\begin{proof}[Proof of \ref{thm5:singularities}] Suppose $S'$ has at least two different balanced completions. These differ by a flow, which in particular implies that the boundary contains a non-branching rational component $F_\8$ with zero self-intersection. Then $F_\8$ is a fiber of a $\PP^1$-ruling $p$ of a balanced completion $(V,D)$. We may assume that $F_\8$ is not contained in any maximal twig of $D$, Indeed, after moving the $0$-curve by a flow to a tip of a new boundary it gives an affine ruling of $S'$, which is possible only if $\ovk(S_0)=-\8$. Since $F_\8$ is non-branching, the induced ruling restricts to an untwisted $\C^*$-ruling of $S'$. It follows from the connectedness of the modification $\eta $ (see the proof of \ref{thm5:(k,k0) computations}) that $n>0$, so this restriction has more than one singular fiber. Moreover, both components of $D_h$ are branching in $D$. Since $F_\8$ is the only non-branching $0$-curve in $D$, centers of elementary transformations lie on the intersection of the fiber at infinity with $D_h$. If $D$ is strongly balanced then one of the components of $D_h$ is a $0$-curve, hence there are at most two strongly balanced completions. Conversely, suppose $S'$ has an untwisted $\C^*$-ruling with base $\C^1$ and $n>0$ and let $(V,D,p)$ be a completion of this ruling. As $S'$ is not affine-ruled, the horizontal components $H$, $H'$ of $D$ are branching, so $(V,D)$ is balanced and we can assume $H'^2=0$. Since $H$, $H'$ are proper transforms of two disjoint sections on a Hirzebruch surface, we have $H^2+H'^2+n\leq 0$, so $H^2\neq 0$ and we can obtain a different strongly balanced completion of $S'$ by a flow which makes $H$ into a $0$-curve.

(2), (3) By \cite[4.5]{Palka-classification1} and \cite{Palka-k(S_0)=0} we may assume that $S'$ is $\C^*$-ruled. If this ruling is untwisted then it follows from the proof of \ref{thm5:(k,k0) computations} that $S'$ has a unique singular point and it is a cyclic singularity. In the twisted case, since $\E\subseteq F_0$, we see that if $\E$ is not connected then $F_0$ is of type (A)(i) and if $\E$ is not a chain then $F_0$ is of type (A)(iv). \end{proof}

\bsrem The set of isomorphism classes of strongly balanced boundaries that a given surface admits is an invariant of the surface, which can easily distinguish between many $\Q$-acyclic surfaces. \esrem

\subsection{Singular planes of negative Kodaira dimension}\label{ssec:S' of negative Kod dim}

As another corollary from \ref{thm5:(k,k0) computations} we give a detailed description of singular $\Q$-homology planes of negative Kodaira dimension. We assume that $\ovk(S_0)\neq 2$, but as we show in \cite{PaKo-general_type} this assumption is in fact redundant.

\bthm\label{thm5:k(S') negative} Let $S'$ be a singular $\Q$-homology plane of negative Kodaira dimension and let $S_0$ be its smooth locus. If $\ovk(S_0)\neq 2$ then exactly one of the following holds:\benum[(i)]

\begin{figure}[h]\centering\includegraphics[scale=0.35]{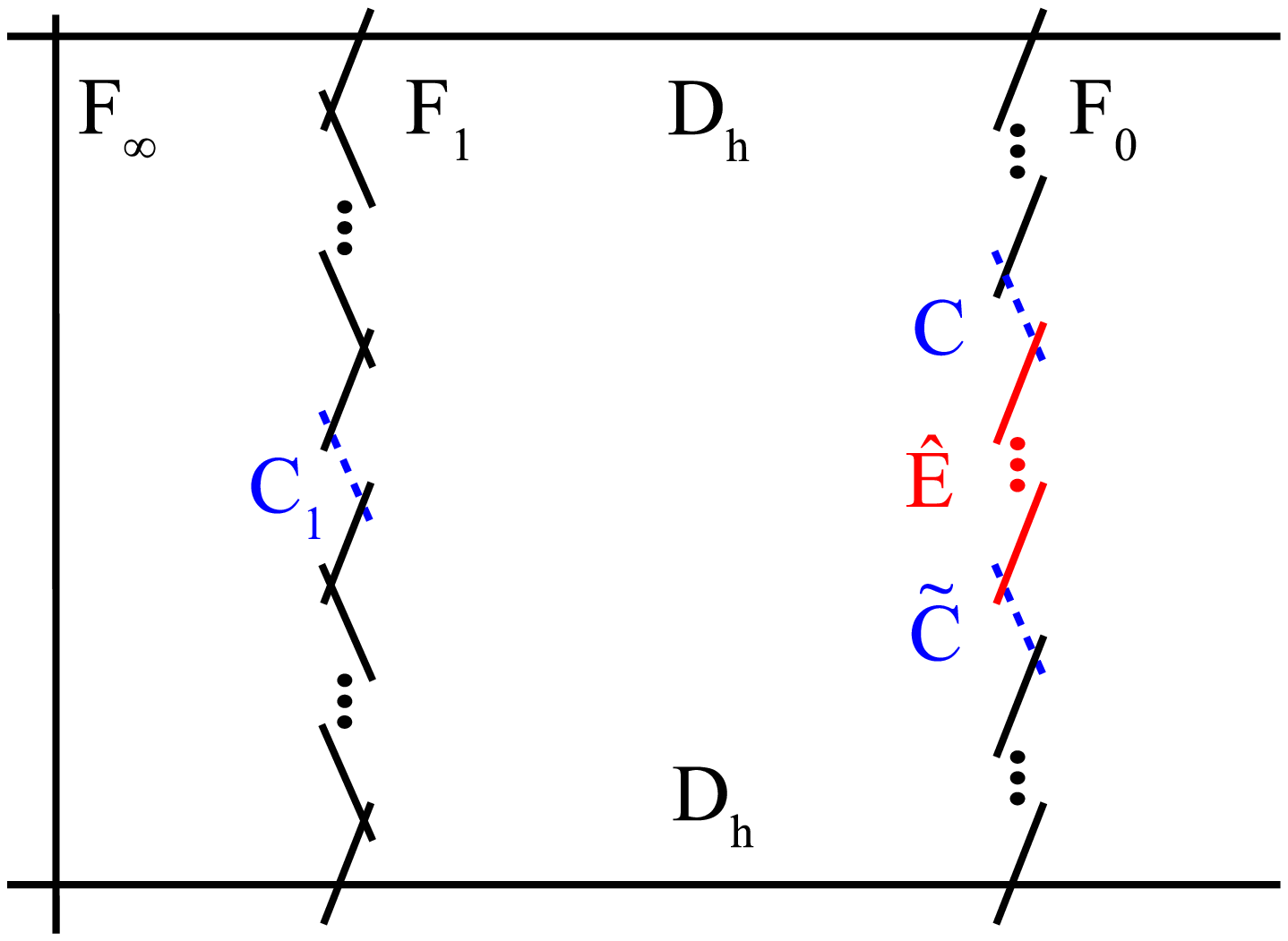}\caption{Untwisted $\C^*$-ruling, $\ovk(S')=-\8$}  \label{fig:extendable}\end{figure}

\item $\ovk(S_0)=-\8$, $S'$ is affine-ruled or isomorphic to $\C^2/G$ for a small finite non-cyclic subgroup $G<GL(2,\C)$,

\item $\ovk(S_0)\in\{0,1\}$, $S'$ is non-logarithmic and is isomorphic to a quotient of an affine cone over a smooth projective curve by an action of a finite group acting freely off the vertex of the cone and preserving the set of lines through the vertex,

\item $\ovk(S_0)\in\{0,1\}$, $S'$ has an untwisted $\C^*$-ruling with base $\C^1$ and two singular fibers, one of them consists of two $\C^1$'s meeting in a cyclic singular point, after taking a resolution and completion the respective completed singular fiber is of type (B)(i) with $\mu,\wt \mu\geq 2$ (see Fig.~\ref{fig:extendable}, cf. \ref{thm5:(k,k0) computations}). \eenum \ethm

\begin{proof} By the results of \cite[4.5]{Palka-classification1}, \cite{Palka-k(S_0)=0} and section \ref{sec:k(S0) negative} we may assume that $S'$ is logarithmic, $\C^*$-ruled and $\ovk(S_0)\geq 0$. We need to show (iii). Let $(V,D,p)$ be a minimal completion of the $\C^*$-ruling. We use \ref{thm5:(k,k0) computations}. If $p$ is twisted then $$0>\kappa_0\geq \lambda-\frac{1}{2}\geq \frac{n-1}{2},$$ so $n=\lambda=0$. The inequalities $\kappa<0$ and $\kappa_0\geq 0$ can be satisfied only in case (A)(iii) and then $D_h^2=0$ by \ref{lem5:minimal model of a Qhp}, so $D_h$ induces an untwisted $\C^*$-ruling of $S'$. Suppose $p$ is untwisted. Since $\kappa\neq \kappa_0$, $p$ has base $\C^1$ and is of type  (B)(i). Since $$0>\kappa=\lambda-1\geq \frac{n}{2}-1,$$ we get $n\leq 1$, but for $n=0$ we get $\kappa_0<\lambda<0$, so in fact $n=1$. Then $0\leq \kappa_0=1-\frac{1}{\mu_1}-\frac{1}{min(\mu,\wt \mu)},$ hence $min(\mu,\wt \mu)\geq 2$. \end{proof}

By \ref{prop:homology} $H_i(S',\Z)$ vanishes for $i>1$. If $S'$ is of type $\C^2/G$ or of type (ii) then it is contractible, $H_1(S',\Z)$ for affine-ruled $S'$ was computed in \ref{rem3:H1 of affine-ruled S'}. For completeness we now compute the fundamental group of $S'$ of type (iii), which by \ref{prop:homology} is the same as $\pi_1(S)$. Let $E_0$ be a component of $\E$ intersecting $C$. Contract $\wt C$ and successive vertical $(-1)$-curves until $C$ is the only $(-1)$-curve in the fiber ($C$ cannot became a $0$-curve, because it does not intersect $D_h$), denote this contraction by $\theta$. Let $\theta'$ be the contraction of $\theta_*F_0$ and $F_1$ to smooth fibers. Put $U=S_0\setminus (C_1\cup C\cup \wt C)$ and let $\gamma_1,\gamma,t\in\pi_1(U)$ be the vanishing loops of the images of $F_1$, $F_0$ under $\theta'\circ\theta$ and of some component of $D_h$  (cf. \cite[4.17]{Fujita}). We need to compute the kernel of the epimorphism $\pi_1(U)\to\pi_1(S)$. Since $\theta$ does not touch $C$, $\theta_*E_0\neq 0$ and $\theta_*F_0$ is columnar. Using 7.17 loc. cit. one can show by induction on the number of components of a columnar fiber that since $E_0\cdot C\neq 0$, the vanishing loops of $E_0$ and $C$, which are of type $\gamma^at^b$ and $\gamma^ct^d$, satisfy $ad-bc=\pm 1$. Thus $\gamma$ and $t$ are in the kernel, hence $$\pi_1(S)=\langle\gamma_1:\gamma^{\mu_1}\rangle\cong \Z_{\mu_1}.$$ In particular, $S'$ is not a $\Z$-homology plane.

\section{Uniqueness of $\C^*$-rulings}\label{ssec:number of C*-rulings}

\subsection{The number of $\C^*$-rulings} We now consider the question of uniqueness of $\C^*$-rulings of $S_0$ and $S'$. Recall that a $\C^*$-ruling of $S_0$ is \emph{extendable} if it extends to a ruling (morphism) of $S'$. Two rational rulings of a given surface are considered the same if they differ by an automorphism of the base. In case a $\C^*$-ruling of $S_0$ exists using the information on snc-minimal boundaries we are able to compute the number of different $\C^*$-rulings.

\bthm\label{thm5:uniqueness of C*-rulings} Let $S'$ be a singular $\Q$-homology plane which is not affine-ruled. Let $p_1,\ldots,p_r$, $r\in \N\cup \{\8\}$ be all different $\C^*$-rulings of the smooth locus $S_0$ of $S'$. Let $D$ be an snc-minimal boundary of $S'$. \benum[(1)]

\item If $\ovk(S_0)=2$ or if $S'$ is exceptional (hence $\ovk(S_0)=0$) then $r=0$.

\item If $\ovk(S_0)=1$ or if $S'$ is non-logarithmic then $r=1$.

\item If $\ovk(S_0)=-\8$ then $r\geq 1$ and $p_1$ is non-extendable. Moreover, $r\neq 1$ only if the fork which is an exceptional divisor of the snc-minimal resolution of $S'$ is of type $(2,2,k)$. In the last case we have:\benum [(i)]

    \item if $k\neq 2$ then $r=2$, $p_2$ is twisted and has a unique singular fiber, which is of type (A)(iv),

    \item if $k=2$ then $r=4$, $p_2,p_3,p_4$ are twisted and they all have unique singular fibers, which are of type (A)(iv). \eenum

\item Assume that $\ovk(S_0)=0$, $S'$ is logarithmic and not exceptional. Then all $p_i$ extend to $\C^*$-rulings of $S'$ and the following hold:\benum[(i)]

   \item If the dual graph of $D$ is $$ \xymatrix{{-2}\ar@{-}[r] &{-1}\ar@{-}[r]\ar@{-}[d]& {k}\ar@{-}[r]\ar@{-}[d] &{-2}\\ {} & {-2} &{-2} & {}} $$ with $k\leq -2$ then $r=1$ and $p_1$ is twisted.

   \item If the dual graph of $D$ is $$ \xymatrix{{-2}\ar@{-}[r] &{-1}\ar@{-}[r]\ar@{-}[d]& {-1}\ar@{-}[r]\ar@{-}[d] &{-2}\\ {} & {-2} &{-2} & {}} $$ then $r=2$ and $p_1$, $p_2$ are twisted.

   \item If the dual graph of $D$ is $$ \xymatrix{{-2}\ar@{-}[r] &{k}\ar@{-}[r]\ar@{-}[d]& {0}\ar@{-}[r] &{m}\ar@{-}[r]\ar@{-}[d] &{-2}\\ {} & {-2} &{} &{-2} & {}} $$ then $r=3$, $p_1$, $p_2$ are twisted and $p_3$ is untwisted with base $\C^1$.

   \item In all other cases $r=2$, $p_1$ is twisted and $p_2$ is untwisted. \eenum\eenum\ethm

\begin{proof} (1) By definition exceptional $\Q$-homology planes are not $\C^*$-ruled. If $S_0$ is of general type then $S_0$ is not $\C^*$-ruled by Iitaka's easy addition formula \cite[10.4]{Iitaka}.

(2) If $S'$ is non-logarithmic then the $\C^*$-ruling of $S'$ is unique by \cite[4.1]{Palka-classification1}. Assume that $\ovk(S_0)=1$. Let $(\ov S,D)$ be some normal completion of the snc-minimal resolution $S\to S'$. Denote the exceptional divisor of the resolution by $\E$. By \cite[6.11]{Fujita} for some $n>0$ the base locus of $|n(K_{\ov S}+D+\E)^+|$ is empty and the linear system gives a $\PP^1$-ruling of $\ov S$ which restricts to a $\C^*$-ruling of $S_0$ (cf. also \cite[2.6.1]{Miyan-OpenSurf}). Consider another $\C^*$-ruling of $S_0$. Modifying $\ov S$ if necessary we can assume that it extends to a $\PP^1$-ruling of $\ov S$. Let $f'$ be a general fiber of this extension. Then $$f'\cdot (K_{\ov S}+D+\E)=f'\cdot K_{\ov S}+2=0,$$ hence $$f'\cdot (K_{\ov S}+D+\E)^+ +f'\cdot (K_{\ov S}+D+\E)^-=0.$$ However, $(K_{\ov S}+D+\E)^-$ is effective and $(K_{\ov S}+D+\E)^+$ is numerically effective, so $$f'\cdot (K_{\ov S}+D+\E)^+=f'\cdot (K_{\ov S}+D+\E)^-=0,$$ and we see that the rulings are the same.

(3), (4) First we need to understand how to find all twisted $\C^*$-rulings of a given $S'$. Consider a twisted $\C^*$-ruling of $S'$ and let $(\wt V,\wt D,\wt p)$ be a minimal completion of this ruling. By the $\wt p$-minimality of $\wt D$, $\wt D_h$ is the only component of $\wt D$ which can be a non-branching $(-1)$-curve, so there is a connected modification $(\wt V,\wt D)\to (V,D)$ with snc-minimal $D$. Let $\wt D_0\subseteq \wt D$ be the $(-1)$-curve of the fiber at infinity (cf. \ref{lem5:F_0 basic properties}). Note that $D$ is not a chain, otherwise $S'$ is affine-ruled. Let $D_0\subseteq D$ be the image of $\wt D_0$ and let $T$ be the connected component of $D-D_0$ containing the image of the horizontal component (which is a point if the modification is nontrivial). In this way a twisted $\C^*$-ruling of $S'$ determines a pair $(D_0,T)$ (with $D_0+T$ contained in a boundary of some snc-minimal completion), such that $\beta_D(D_0)=3$, $D_0^2\geq -1$, $T$ is a connected component of $D-D_0$ containing the image of the horizontal section and both connected components of $D-D_0-T$ are $(-2)$-curves. Conversely, if we have an snc-minimal normal completion $(V,D)$ and a pair as above, we make a connected modification $(\wt V,\wt D)\to (V,D)$ over $D$ by blowing successively on the intersection of the total transform of $T$ with the proper transform of $D_0$ until $D_0$ becomes a $(-1)$-curve. The $(-1)$-curve together with the transform of $D-T-D_0$ induce a $\PP^1$-ruling of $V'$ and constitute the fiber at infinity for this ruling. The restriction to $S'$ is a twisted $\C^*$-ruling.

Suppose $\ovk(S_0)=-\8$. Since $S_0$ is not affine-ruled, $S'\cong \C^2/G$ for a finite noncyclic small subgroup $G<GL(2,\C)$ (cf. section \ref{sec:k(S0) negative}). Let $(V,D)$ be an snc-minimal normal completion of $S'$ and let $\ov S\to V$ be a minimal resolution with exceptional divisor $\E$. We saw in the proof of \ref{prop4:non affine-ruled S'} that $S_0$ admits a Platonic $\C^*$-ruling, which extends to a $\PP^1$-ruling of $\ov S$. Moreover, $D$ and $\E$ are forks for which $D_h$ and $\E_h$ are the unique branching components of $D$ and $E$ respectively. In particular, the $\C^*$-ruling does not extend to a ruling of $S'$ and as non-branching components of $D$ have negative self-intersections, $(\ov S,D+\E)$ is a unique snc-minimal smooth completion of $S_0$ (and hence $(V,D)$ is a unique snc-minimal normal completion of $S'$). It follows from the proof of \cite[4.1]{Palka-classification1} that the non-extendable $\C^*$-ruling of $S_0$ is unique. Suppose there is a $\C^*$-ruling of $S_0$ which does extend to $S'$. Since $\E$ is not a chain, it follows from the proof of \ref{thm5:(k,k0) computations} that this ruling is twisted. Since maximal twigs of $\E$ and $D$ are adjoint chains of columnar fibers, we see that a maximal twig of $D-D_h$ is a $(-2)$-curve if and only if the respective maximal twig of $\E-\E_h$ is a $(-2)$-curve. Moreover, we have $0<d(\E)$, so $\E_h^2\leq -2$ and since $\E_h^2+D_h^2=-3$, we have $D_h^2\geq -1$. Therefore, $S'$ admits a twisted $\C^*$-ruling if and only if $\E$ is a fork of type $(2,2,k)$ for some $k\geq 2$. If $k\neq 2$ then the choice of $(D_0,T)$ as above is unique and if $k=2$ then there are three such choices. Note that if $(V',D',p)$ is a minimal completion of such a ruling then $D'$ is a fork, so since $\kappa_0<0$, we have $n=0$ and $F_0$ is of type (A)(iv) (cf. the proof of \ref{thm5:(k,k0) computations}). This gives (3).

We can now assume that $\ovk(S_0)=0$, $S'$ is logarithmic and not exceptional. Then $S_0$ is $\C^*$-ruled and by \cite[4.7(iii)]{Palka-classification1} each $\C^*$-ruling of $S_0$ extends to a $\C^*$-ruling of $S'$. Let $r\in\{1,2,\ldots\}\cup \{\8\}$ be the number of all different (up to automorphism of the base) $\C^*$-rulings of $S'$ and let $(V_i,D_i,p_i)$ for $i\leq r$ be their minimal completions. Minimality implies that non-branching $(-1)$-curves in $D_i$ are $p_i$-horizontal. We add consequently an upper index $(i)$ to objects defined previously for any $\C^*$-ruling when we refer to the ruling $p_i$. If $p_i$ is untwisted we denote the horizontal components of $D_h^{(i)}$ by $H^{(i)}$, $H'^{(i)}$.

Suppose $p_1$ is untwisted with base $\PP^1$. Then $F_0^{(1)}$ contains a rivet and by \ref{cor5:when k_0=0} $n^{(1)}=2$, so $D_1$ does not contain non-branching $b$-curves with $b\geq -1$. Then $(V_1,D_1)$ is balanced and $S'$ does not admit an untwisted $\C^*$-ruling with base $\C^1$, as it does not contain non-branching $0$-curves (cf. \ref{lem5:F_0 basic properties}). By \ref{cor5:when k_0=0} each component of $D_h^{(1)}$ has $\beta_{D_1}=3$ and intersects two $(-2)$-tips of $D_1$. Note that $\zeta^{(1)}$ (cf. \ref{lem5:minimal model of a Qhp}) touches $D_h^{(1)}$ two times if both components of $D_h^{(1)}$ intersect the same horizontal component of $F_0^{(1)}$ and three times if not. By \ref{lem5:minimal model of a Qhp} and by the properties of Hirzebruch surfaces we get $-3\leq (D_h^{(1)})^2\leq -2$. In particular, one of the components of $D_h^{(1)}$, say $H^{(1)}$, has $(H^{(1)})^2\geq -1$, so by the discussion about twisted $\C^*$-rulings above $H^{(1)}$ together with two $(-2)$-tips of $D_1$ gives rise to a twisted $\C^*$-ruling $p_2$ of $S'$. Since $H'^{(1)}$ together with two $(-2)$-tips of $D_1$ intersecting it are contained in a fiber of $p_2$, $(H'^{(1)})^2\leq -2$. Thus $p_2$ is the only twisted ruling of $S'$, because $H^{(1)}$ is the only possible choice for a middle component of the fiber at infinity of a twisted ruling. Suppose $r\geq 3$. Then $p_3$ is untwisted with base $\PP^1$. Since $D_1$ does not contain non-branching $0$-curves, any flow in $D_1$ is trivial, so $V_3=V_1$. Since $p_3$ and $p_1$ are different after restriction to $S'$, the $S_0$-components $C^{(1)}$, $C^{(3)}$ contained respectively in $F_0^{(1)}$, $F_0^{(3)}$ are different. As they both intersect $\E$, they are contained in the same fiber of $p_2$, a contradiction with $\Sigma_{S_0}^{(2)}=0$. Note that since $D$ contains no non-branching $0$-curves, $D$ is not of type (iii). Since $n^{(1)}=2$, $D$ contains at least seven components, so $D$ is not of type (i) or (ii).

We can now assume that each untwisted $\C^*$-ruling of $S'$ has base $\C^1$. Suppose $p_1$ is such a ruling. By \ref{cor5:when k_0=0} both horizontal components of $D_1$ have $\beta_{D_1}=3$ and one of them, say $H'^{(1)}$, intersects two $(-2)$-tips $T$ and $T'$ of $D_1$. In particular, $D_1$ is snc-minimal. Since $F_\8^{(1)}=[0]$, changing $V_1$ by a flow if necessary we may assume that $H'^{(1)}$ is a $(-1)$-curve. Then $F_\8^{(2)}=T+2H'^{(1)}+T'$ induces a $\PP^1$-ruling $p_2:V_1\to \PP^1$, which is a twisted $\C^*$-ruling after restricting it to $S'$. Suppose $r\geq 3$. If $p_3$ is untwisted then its base is $\C^1$ and changing $V_3$ by a flow if necessary we can assume that $V_3=V_1$. But then $F_\8^{(1)}=F_\8^{(3)}$, because $D_1$ contains only one non-branching $0$-curve, so $p_1$ and $p_3$ have a common fiber and hence cannot be different after restriction to $S'$, a contradiction. Thus $p_3$ is twisted. By the discussion above $p_3$ can be recovered from a pair $(D_0,T)$ on some snc-minimal completion of $S'$. All such completions of $S'$ differ from $(V_1,D_1)$ by a flow, which is an identity on $V_1-F_\8^{(1)}$, hence the birational transform of $D_0$ on $V_1$ is either $H^{(1)}$ or $H'^{(1)}$. Since the restrictions of $p_1$ and $p_2$ to $S'$ are different, it is $H^{(1)}$. It follows that $r=3$ and $D_1-H'^{(1)}$ has two $(-2)$-tips as connected components, hence the dual graph of $D_1$ is as in (iii). Conversely, if $S'$ has a boundary as in (iii) then besides the untwisted $\C^*$-ruling induced by the $0$-curve it has also two twisted rulings, each with one of the branching components as the middle component of the fiber at infinity.

We can finally assume that all $\C^*$-rulings of $S'$ are twisted. Let $(V,D)$ be a balanced completion of $S'$. Since $S'$ does not admit untwisted $\C^*$-rulings, $D$ does not contain non-branching $0$-curves, so $(V,D)$ is a unique snc-minimal completion of $S'$. Thus to find all twisted $\C^*$-rulings of $S'$ we need to determine all pairs $(D_0,T)$, such that $D_0+T\subseteq D$,  $D_0^2\geq -1$, $\beta_D(D_0)=3$ and $D-T-D_0$ consists of two $(-2)$-tips. Let $(D_0,T)$ and $(D_0',T')$ be two such pairs. Suppose $D_0\neq D_0'$ and, say, $D_0'^2\geq D_0^2$. We have $D_0\cdot D_0'\neq 0$, otherwise the chain $D-T'$, which is not negative definite, would be contained (and not equal, since $\nu\leq 1$) in a fiber of the twisted ruling associated with $(D_0,T)$, which is impossible. Then $D$ has six components and we check that $$d(D)=16((D_0^2+1)(D_0'^2+1)-1),$$ so $(D_0^2+1)(D_0'^2+1)\leq 0$, because $d(D)<0$. Then $D_0^2=-1$ and $D_0'$ is a 2-section of the twisted ruling associated with $(D_0,T)$. Since $\beta_D(D_0')=3$, by \ref{cor5:when k_0=0} and \ref{lem5:minimal model of a Qhp} for this ruling $n=1$, $D_0'$ is a $(-1)$-curve and $D$ has dual graph as in (ii). Conversely, it is easy to see that $S'$ with such a boundary has two twisted $\C^*$-rulings. Therefore, we can assume that the choice of $D_0$ for a pair $(D_0,T)$ as above is unique. Let $p_1$ be a twisted $\C^*$-ruling associated with some pair $(D_0,T)$. Suppose $n^{(1)}=0$. By \ref{lem5:minimal model of a Qhp} $\zeta_*D_h^{(1)}$ is a $0$-curve, so $F=\zeta^*\zeta_*D_h^{(1)}$ induces a $\PP^1$-ruling $p$ of $V$. If $\zeta$ touches $D_h^{(1)}$ then $F$ contains the $S_0$-component of $F_0^{(1)}$, so $F\nsubseteq D$ and $p$ restricts to an untwisted $\C^*$-ruling of $S'$ with base $\PP^1$. If $\zeta$ does not touch $D_h^{(1)}$ then $p$ restricts to a $\C^*$-ruling of $S'$ with base $\C^1$. This contradicts the assumption. By \ref{cor5:when k_0=0} we get that $n^{(1)}=1$, $F_0^{(1)}$ contains no $D_1$-components and $\mu_1=2$. In particular, $D_1=D$. Moreover, as $n^{(1)}=1$, by \ref{lem5:minimal model of a Qhp} $(D_h^{(1)})^2\leq -1$, so $D$ has a dual graph as in (i) or (ii). Conversely, if $D$ is of type (i) or (ii) then $r=2$ if $k=-1$ and $r=1$ if $k\leq -2$. \end{proof}

\subsection{The number of affine lines}

The theorem \ref{thm5:uniqueness of C*-rulings} has interesting consequences. Namely, it is known (\cite{Zaid_isotrivial}, \cite{GM-Affine-lines}) that $\Q$-homology planes with smooth locus of general type (in particular the smooth ones) do not contain topologically contractible curves. In fact the number $\ell\in \N\cup\{\8\}$ of contractible curves on a $\Q$-homology plane $S'$ is known except two cases: when $S'$ is non-logarithmic or when $S'$ is singular and $\ovk(S_0)=0$ (cf. \cite[10.1]{Palka-recent_progress} and references there). Clearly, in the first case $\ell=\8$ by the main result of \cite{Palka-classification1}. The case when $S'$ is smooth and of Kodaira dimension zero has been considered in \cite{GP-affine_lines_k=0}. The theorem \ref{cor5:number of contractible curves} is the missing piece of information. The methods can be easily applied to recover the result in loc. cit.

\begin{proof}[Proof of \ref{cor5:number of contractible curves}] We can assume that $S'$ is logarithmic. Suppose $S'$ contains a topologically contractible curve $L$. We show that $L$ is vertical for some $\C^*$-ruling of $S'$. Note that the proper transform of $L$ on $\ov S$ meets each connected component of $\E$ in at most one point. We use the logarithmic Bogomolov-Miyaoka-Yau inequality as in \cite[2.12]{KR-ContrSurf} to show that $\ovk(S_0-L)\leq 1$. In case $\ovk(S_0-L)=1$ the surface $S_0-L$ is $\C^*$-ruled (cf. \cite[6.11]{Fujita}), so we may assume that $\ovk(S_0-L)=0$. Let $\Z[D+\E]$ be a free abelian group generated by the components of $D+\E$. Since $$\Pic S_0 = \Coker(\Z[D+\E]\to \Pic \ov S)$$ is torsion, the class of $L$ in $\Pic S_0$ is torsion. Then there exists a morphism $f\:S_0-L\to\C^*$ and taking its Stein factorization we get a $\C^*$-ruling of $S_0-L$, which (as $\ovk(S_0)\neq -\8$) extends to a $\C^*$-ruling of $S_0$. Since $S_0$ is logarithmic, each $\C^*$-ruling of $S_0$ extends in turn to a $\C^*$-ruling of $S'$. Therefore $L$ is vertical for some $\C^*$-ruling of $S'$ and we are done. In particular, exceptional $\Q$-homology planes do not contain contractible curves. It follows from \ref{cor5:when k_0=0} that if the ruling is twisted or untwisted with base $\PP^1$ then the vertical contractible curve is unique and is contained in the unique singular non-columnar fiber. For an untwisted ruling with base $\C^1$ there are at most two such curves. In particular, in cases (4)(i) and (ii) of the theorem \ref{thm5:uniqueness of C*-rulings} $L$ needs to intersect the horizontal component of the boundary, so we get respectively $\ell=1$ and $\ell=2$. In case (4)(iii) the unique vertical contractible curves for the twisted rulings $p_1$ and $p_3$ are distinct and do not intersect the horizontal components of respective rulings, hence are both vertical for the untwisted ruling $p_3$, so $\ell=2$. In the remaining case $(4)(iv)$ we have $r=2$, $p_1$ is twisted and $p_2$ is untwisted. We can assume that the base of $p_2$ is $\C^1$ and the unique non-columnar singular fiber contains two contractible curves, $L_1$ and $L_2$, otherwise $\ell\leq 2$ from the above remarks and we are done. Since the twisted ruling is unique, there is exactly one horizontal component $H$ of $D_h^{(2)}$ which meets two $(-2)$-tips of $D_h^{(1)}$ (together with these tips it induces the twisted ruling). Clearly, only one $L_i$ can intersect $H$, so the second one is vertical for $p_1$ and we get $\ell\leq 2$ is this case too. \end{proof}

\bibliographystyle{amsalpha}
\bibliography{bibl}

\end{document}